\documentclass[11pt,leqno]{article}
\normalfont
\renewcommand{\baselinestretch}{1.1455}
\usepackage{amsfonts}

\newfont{\eulercursive}{eurm10 at 11pt}
\newcommand{\myl}{\mbox{\eulercursive `}}

\newcommand{\QED}{\raisebox{0.5mm}{\fbox{\rule{0mm}{1.5mm}\ }}}

\newcounter{myfn}[page]

\newcommand{\myfootnote}[1]{\setcounter{footnote}{\value{myfn}}%
    \footnote{#1}\stepcounter{myfn}}

\newcommand{\ECoxeterGraphFigure}{Figure 2.1} 
\newcommand{\BtwoMarsDemo}{Figure 2.2} 

\newcommand{\StrongConvergenceTheorem}{Theorem 2.1}
\newcommand{\StrongConvergenceCorollary}{Lemma 2.2}
\newcommand{\ComparisonTheorem}{Theorem 2.3}
\newcommand{\ComparisonCorollary}{Lemma 2.4}
\newcommand{\ComparisonResults}{Lemmas 2.2 and 2.4}
\newcommand{\NotMarsFriendlyLemma}{Lemma 2.5}
\newcommand{\EveryNodeFiredLemma}{Lemma 2.6}
\newcommand{\NewLemmaList}{Lemmas 2.2 and 2.6}
\newcommand{\SubgraphLemma}{Lemma 2.7}
\newcommand{\ErikssonWordProposition}{Theorem 2.8}
\newcommand{\HumphreysTheorem}{Proposition 2.9}
\newcommand{\TFAE}{Proposition 2.10}
\newcommand{\PropList}{Propositions 2.9 and 2.10}
\newcommand{\PositiveToNegativeRootsLemma}{Lemma 2.11}
\newcommand{\LengthProposition}{Proposition 2.12}
\newcommand{\EquivalenceRemark}{Remark 2.13}
\newcommand{\HRTResult}{Proposition 2.14} 
\newcommand{\TitsConeConvergenceResult}{Theorem 2.15}
\newcommand{\RuleOutInfinitePropositions}{Proposition 2.14 and 
Theorem 2.15}

\newcommand{\CoxeterTwoGeneratorAnalysis}{Lemma 3.1}
\newcommand{\CoxeterHumphreysTheorem}{Theorem 3.3}
\newcommand{\CoxeterTFAE}{Theorem 3.6}
\newcommand{\CoxeterPositiveToNegativeRootsLemma}{Lemma 3.8}
\newcommand{\CoxeterLengthProposition}{Theorem 3.9}
\newcommand{\CoxeterHRTResult}{Theorem 4.3} 

\newcommand{\KeyWJResultLemma}{Lemma 3.1}
\newcommand{\KeyWJResult}{Proposition 3.2}
\newcommand{\LengthUsingLongestWord}{Corollary 3.3} 
\newcommand{\KeyWJResultCorollary}{Corollary 3.4}

\newcommand{\AdjacencyFreeForFinite}{Proposition 4.2}
\newcommand{\AdjacencyFreeClassification}{Theorem 4.3}

\newcommand{\NoRepeatPosRootsProp}{Lemma 5.1}

\newcommand{\AllPositiveRootsProp}{Theorem 5.2}

\newcommand{\EGCMTheorem}{Theorem 6.1} 
\newcommand{\ErikssonTheorem}{Theorem 6.2} 
\newcommand{\TwoTheorems}{Theorems 6.1 and 6.2} 
\newcommand{\DeodharEquivalenceCorollary}{Corollary 6.3}
\newcommand{\CycleLemma}{Lemma 6.4}
\newcommand{\FourCycleLemma}{Lemma 6.5}
\newcommand{\ThreeCycleLemma}{Lemma 6.6}
\newcommand{\RuleOutSmallCyclesLemmas}{Lemmas 6.4 and 6.5}

\newcommand{\IntroNum}{1}
\newcommand{\PrelimNum}{2}
\newcommand{\ConvergenceNum}{3}
\newcommand{\AdjacencyFreeNum}{4}
\newcommand{\RootsNum}{5}
\newcommand{\DynkinClassNum}{6}

\newcommand{\selt}{\mathbf{s}} \newcommand{\telt}{\mathbf{t}}

\newcommand{\disjointunion}{\setlength{\unitlength}{0.14cm}
\ 
\begin{picture}(2,2) 
\put(0,0){$\cup$}
\put(0.9675,1.5){\circle*{0.5}}
\end{picture}\ }

\newcommand{\CircleInteger}[1]{
\setlength{\unitlength}{0.14cm}
\begin{picture}(3,2) 
\put(1,1){\circle{2}}
\put(0.55,0.6){{\tiny #1}}
\end{picture}
}

\newcommand{\CircleIntegerm}{
\setlength{\unitlength}{0.14cm}
\begin{picture}(3,2) 
\put(1,1){\circle{2}}
\put(0.3,0.6){{\tiny $m$}}
\end{picture}
}

\newcommand{\AnEGraph}{\setlength{\unitlength}{0.75in}
\begin{picture}(6,0.85)
\put(0,0){\begin{picture}(1,0)
            \put(0,0.1){\circle*{0.075}}
            \put(1,0.1){\circle*{0.075}}
            \put(2,0.1){\circle*{0.075}}
            \put(4,0.1){\circle*{0.075}}
            \put(5,0.1){\circle*{0.075}}
            \put(6,0.1){\circle*{0.075}}
            \put(0,0.1){\line(1,0){2}}
            \multiput(2,0.1)(0.4,0){5}{\line(1,0){0.2}}
            \put(4,0.1){\line(1,0){2}}
           \end{picture}}
\end{picture}}

\newcommand{\BnEGraph}{\setlength{\unitlength}{0.75in}
\begin{picture}(6,0.55)
\put(0,0){\begin{picture}(1,0)
            \put(0,0.1){\circle*{0.075}}
            \put(1,0.1){\circle*{0.075}}
            \put(2,0.1){\circle*{0.075}}
            \put(4,0.1){\circle*{0.075}}
            \put(5,0.1){\circle*{0.075}}
            \put(6,0.1){\circle*{0.075}}
            \put(0,0.1){\line(1,0){2}}
            \multiput(2,0.1)(0.4,0){5}{\line(1,0){0.2}}
            \put(4,0.1){\line(1,0){2}}
            \put(5.3,0.15){\CircleInteger{4}}
           \end{picture}}
\end{picture}}

\newcommand{\DnEGraph}{\setlength{\unitlength}{0.75in}
\begin{picture}(6,0.75)
\put(0,-0.25){\begin{picture}(1,0)
            \put(0,0.35){\circle*{0.075}}
            \put(1,0.35){\circle*{0.075}}
            \put(2,0.35){\circle*{0.075}}
            \put(4,0.35){\circle*{0.075}}
            \put(5,0.35){\circle*{0.075}}
            \put(6,0.1){\circle*{0.075}}
            \put(6,0.6){\circle*{0.075}}
            \put(0,0.35){\line(1,0){2}}
            \multiput(2,0.35)(0.4,0){5}{\line(1,0){0.2}}
            \put(4,0.35){\line(1,0){1}}
            \put(5,0.35){\line(4,1){1}}
            \put(5,0.35){\line(4,-1){1}}
           \end{picture}}
\end{picture}}

\newcommand{\EEightEGraph}{\setlength{\unitlength}{0.75in}
\begin{picture}(6,0.75)
\put(0,-0.25){\begin{picture}(1,0)
            \put(0,0.1){\circle*{0.075}}
            \put(1,0.1){\circle*{0.075}}
            \put(2,0.1){\circle*{0.075}}
            \put(2,0.6){\circle*{0.075}}
            \put(3,0.1){\circle*{0.075}}
            \put(4,0.1){\circle*{0.075}}
            \put(5,0.1){\circle*{0.075}}
            \put(6,0.1){\circle*{0.075}}
            \put(0,0.1){\line(1,0){6}}
            \put(2,0.1){\line(0,1){0.5}}
           \end{picture}}
\end{picture}}

\newcommand{\ESevenEGraph}{\setlength{\unitlength}{0.75in}
\begin{picture}(6,0.75)
\put(0,-0.25){\begin{picture}(1,0)
            \put(0,0.1){\circle*{0.075}}
            \put(1,0.1){\circle*{0.075}}
            \put(2,0.1){\circle*{0.075}}
            \put(2,0.6){\circle*{0.075}}
            \put(3,0.1){\circle*{0.075}}
            \put(4,0.1){\circle*{0.075}}
            \put(5,0.1){\circle*{0.075}}
            \put(0,0.1){\line(1,0){5}}
            \put(2,0.1){\line(0,1){0.5}}
            \put(4.95,0.2){\small {\tt *}}
           \end{picture}}
\end{picture}}

\newcommand{\ESixEGraph}{\setlength{\unitlength}{0.75in}
\begin{picture}(6,0.75)
\put(0,-0.25){\begin{picture}(1,0)
            \put(0,0.1){\circle*{0.075}}
            \put(1,0.1){\circle*{0.075}}
            \put(2,0.1){\circle*{0.075}}
            \put(2,0.6){\circle*{0.075}}
            \put(3,0.1){\circle*{0.075}}
            \put(4,0.1){\circle*{0.075}}
            \put(0,0.1){\line(1,0){4}}
            \put(2,0.1){\line(0,1){0.5}}
            \put(-0.05,0.2){\small {\tt *}}
            \put(3.95,0.2){\small {\tt *}}
           \end{picture}}
\end{picture}}

\newcommand{\FFourEGraph}{\setlength{\unitlength}{0.75in}
\begin{picture}(1,0.85)
\put(0,0){\begin{picture}(1,0)
            \put(0,0.1){\circle*{0.075}}
            \put(1,0.1){\circle*{0.075}}
            \put(2,0.1){\circle*{0.075}}
            \put(3,0.1){\circle*{0.075}}
            \put(0,0.1){\line(1,0){3}}
            \put(1.3,0.15){\CircleInteger{4}}
            \end{picture}}
\end{picture}}

\newcommand{\HFourEGraph}{\setlength{\unitlength}{0.75in}
\begin{picture}(1,0.55)
\put(0,0){\begin{picture}(1,0)
            \put(0,0.1){\circle*{0.075}}
            \put(1,0.1){\circle*{0.075}}
            \put(2,0.1){\circle*{0.075}}
            \put(3,0.1){\circle*{0.075}}
            \put(0,0.1){\line(1,0){3}}
            \put(0.3,0.15){\CircleInteger{5}}
            \end{picture}}
\end{picture}}

\newcommand{\HThreeEGraph}{\setlength{\unitlength}{0.75in}
\begin{picture}(1,0.55)
\put(0,0){\begin{picture}(1,0)
            \put(0,0.1){\circle*{0.075}}
            \put(1,0.1){\circle*{0.075}}
            \put(2,0.1){\circle*{0.075}}
            \put(0,0.1){\line(1,0){2}}
            \put(0.3,0.15){\CircleInteger{5}}
            \put(1.95,0.2){\small {\tt *}}
            \end{picture}}
\end{picture}}

\newcommand{\ITwoEGraph}{\setlength{\unitlength}{0.75in}
\begin{picture}(1,0.55)
\put(0,0){\begin{picture}(1,0)
            \put(0,0.1){\circle*{0.075}}
            \put(1,0.1){\circle*{0.075}}
            \put(0,0.1){\line(1,0){1}}
            \put(0.3,0.15){\CircleIntegerm}
            \end{picture}}
\end{picture}}

\newcommand{\BTwoGraphForFigure}{\setlength{\unitlength}{0.6in}
\begin{picture}(1,0.4)
\put(0,0){\begin{picture}(1,0)
            \put(0,0.1){\circle*{0.075}}
            \put(-0.11,-0.05){\scriptsize $\gamma_{1}$}
            \put(1,0.1){\circle*{0.075}}
            \put(1,-0.05){\scriptsize $\gamma_{2}$}
            \put(0,0.1){\line(1,0){1}}
            \put(0.2,0.1){\vector(1,0){0.1}}
            \put(0.8,0.1){\vector(-1,0){0.1}}
            \put(0.7,0.1){\vector(-1,0){0.1}}
            \end{picture}}
\end{picture}}

\newcommand{\BThreeGraph}{\setlength{\unitlength}{0.75in}
\begin{picture}(2.25,0.275)
\put(0.05,0.05){\begin{picture}(1,0)
            \put(0,0.1){\circle*{0.05}}
            \put(-0.05,-0.05){\footnotesize $\gamma_{1}$}
            \put(0,0.1){\line(1,0){1}}
            \end{picture}}
\put(1.05,0.05){\begin{picture}(1,0)
            \put(0,0.1){\circle*{0.05}}
            \put(-0.05,-0.05){\footnotesize $\gamma_{2}$}
            \put(1,0.1){\circle*{0.05}}
            \put(0.95,-0.05){\footnotesize $\gamma_{3}$}
            \put(0,0.1){\line(1,0){1}}
            \put(0.2,0.1){\vector(1,0){0.1}}
            \put(0.235,-0.025){\tiny 1}
            \put(0.8,0.1){\vector(-1,0){0.1}}
            \put(0.71,-0.025){\tiny 2}
            \end{picture}}
\end{picture}}

\newcommand{\TwoCitiesGraphWithLabels}{
\setlength{\unitlength}{1in}
\begin{picture}(1,0.3)
\put(0,0.15){\begin{picture}(1,0)
            \put(0,0.1){\circle*{0.05}}
            \put(-0.15,-0.05){\large $\gamma_{1}$}
            \put(1,0.1){\circle*{0.05}}
            \put(1.05,-0.05){\large $\gamma_{2}$}
            \put(0,0.1){\line(1,0){1}}
            \put(0.2,0.1){\vector(1,0){0.1}}
            \put(0.8,0.1){\vector(-1,0){0.1}}
            \put(0.225,0){\footnotesize $p$}
            \put(0.71,0){\footnotesize $q$}
            \end{picture}}
\end{picture}}

\newcommand{\WeirdATwoEGCMGraphh}{
\setlength{\unitlength}{1in}
\begin{picture}(1.25,0.195)
\put(0,0.05){\begin{picture}(1,0)
            \put(0,0.1){\circle*{0.05}}
            \put(-0.15,-0.05){\large $\gamma_{1}$}
            \put(1,0.1){\circle*{0.05}}
            \put(1.05,-0.05){\large $\gamma_{2}$}
            \put(0,0.1){\line(1,0){1}}
            \put(0.2,0.1){\vector(1,0){0.1}}
            \put(0.8,0.1){\vector(-1,0){0.1}}
            \put(0.225,-0.025){\scriptsize $1/2$}
            \put(0.71,-0.025){\scriptsize $2$}
            \end{picture}}
\end{picture}}

\newcommand{\ATwoGraphNoEdgeLabels}{\setlength{\unitlength}{0.75in}
\begin{picture}(1.65,0.25)
\put(0.25,0){\begin{picture}(1,0)
            \put(0,0.1){\circle*{0.05}}
            \put(-0.20,-0.05){\large $\gamma_{1}$}
            \put(1,0.1){\circle*{0.05}}
            \put(1.05,-0.05){\large $\gamma_{2}$}
            \put(0,0.1){\line(1,0){1}}
            \end{picture}}
\end{picture}}

\newcommand{\EGCMGraphCirclem}{\setlength{\unitlength}{0.75in}
\begin{picture}(1.65,0.25)
\put(0.25,0){\begin{picture}(1,0)
            \put(0,0.1){\circle*{0.05}}
            \put(-0.20,-0.05){\large $\gamma_{1}$}
            \put(1,0.1){\circle*{0.05}}
            \put(1.05,-0.05){\large $\gamma_{2}$}
            \put(0,0.1){\line(1,0){1}}
            \put(0.3,0.15){\CircleIntegerm}
            \end{picture}}
\end{picture}}

\begin{document}

\newpage
\setcounter{page}{1} 
\renewcommand{\baselinestretch}{1}

\begin{center}
{\large \bf Eriksson's numbers game and finite Coxeter groups}

Robert G.\ Donnelly

Department of Mathematics and Statistics, Murray State
University, Murray, KY 42071

{\bf To appear in the {\em European Journal of Combinatorics}}
\end{center}

\vspace*{-0.25in}
\begin{abstract}
The numbers game is a one-player game played on  
a finite simple 
graph with certain ``amplitudes'' assigned to its edges 
and with an initial assignment of real 
numbers to its nodes.  The moves of the game successively transform 
the numbers at the nodes using the amplitudes in a certain way.  
This game and its interactions with Coxeter/Weyl group theory 
and Lie theory have been studied by many authors.  
In particular, Eriksson connects certain geometric representations of 
Coxeter groups with games on graphs with certain real number 
amplitudes. 
Games played on such graphs are ``E-games.'' 
Here we investigate various finiteness aspects of E-game play:   
We extend Eriksson's work relating moves of the game to reduced 
decompositions of elements of a Coxeter group naturally associated to 
the game graph.    
We use 
Stembridge's theory of fully commutative Coxeter group elements to 
classify what we call here the ``adjacency-free'' initial positions 
for finite E-games.  
We characterize when 
the positive roots for certain geometric representations 
of finite Coxeter groups can be obtained from E-game play. 
Finally, we provide a new Dynkin diagram classification result 
of E-game graphs  
meeting a certain finiteness requirement. 
\begin{center}

\ 
\vspace*{-0.1in}


\end{center}
\end{abstract}

%==================================================================

\noindent
{\Large \bf \IntroNum.\ \ Introduction}

\vspace{1ex} 
The numbers game is a one-player game played on a finite simple graph 
with weights (which we call ``amplitudes'') on its edges 
and with an initial assignment of 
real numbers 
to its nodes.  
Each of the two edge amplitudes (one for each direction) 
will be certain negative real numbers. 
The 
move a player can make 
is to ``fire'' one of the nodes with a positive number.  This move 
transforms the number at the fired node 
by changing its sign, and it also 
transforms the number at each adjacent node in a certain way 
using an amplitude   
along the incident edge.  
The player fires the nodes in some sequence of 
the player's choosing, continuing until no node has a positive 
number.  

The numbers game has been an object of interest for many authors. 
For graphs with integer 
amplitudes the game is attributed to Mozes \cite{Mozes}.  
Eriksson has studied the game 
extensively, see for example \cite{ErikssonLinear}, \cite{ErikssonThesis}, 
\cite{ErikssonJerusalem}, \cite{ErikssonReachable}, 
\cite{ErikssonDiscrete}, \cite{ErikssonEur}, \cite{DE}. 
Eriksson's numbers game allows for certain real number amplitudes.  
Particularly 
important for this paper is his ground-breaking work in 
\cite{ErikssonThesis}, 
\cite{ErikssonDiscrete}, and \cite{ErikssonEur} analyzing  
convergence of numbers games and of 
the connection between the numbers game and Coxeter groups.  
Much of the numbers game discussion in \S 4.3 of the book \cite{BB} by 
Bj\"{o}rner and Brenti can be found in  
\cite{ErikssonThesis} and \cite{ErikssonDiscrete}.  
The game has also been 
studied by Proctor 
\cite{PrEur}, \cite{PrDComplete}, Bj\"{o}rner 
\cite{Bjorner}, and 
Wildberger 
\cite{WildbergerAdv}, \cite{WildbergerEur}, \cite{WildbergerPreprint}.  
Wildberger studies a dual version which 
he calls the ``mutation game.''  See Alon {\em et al} \cite{AKP} for a 
brief and readable treatment of the numbers game on ``unweighted'' 
cyclic graphs. 
The numbers game 
facilitates computations with Coxeter groups and their geometric 
representations (e.g.\ see \S 4.3 of \cite{BB} or \S 3, 4 below).  
Proctor developed this process in \cite{PrEur} 
to compute Weyl group 
orbits of weights with respect to the fundamental weight basis.  Here 
we use his perspective of firing 
nodes with positive, as opposed to negative, numbers. 
In \cite{DW}, we use data from certain numbers games to obtain 
distributive lattice models for families of semisimple Lie algebra 
representations and their Weyl characters.  

This paper extends Eriksson's work, focussing on play from ``dominant'' 
positions where all numbers are either fireable or zero, for which 
the connection to Coxeter groups turns out to be quite explicit.  We 
will let $J$ denote the set  
of nodes where the numbers are zero, and $J^{c}$ denotes its complement. 
The main results can be summarized as follows: 
In \S \ConvergenceNum\ we show how, under 
certain finiteness assumptions, legal play sequences 
from a $J^{c}$-dominant position correspond to reduced words in the 
quotient $W^{J}$. 
In \S \AdjacencyFreeNum\ we relate Stembridge's notion of full 
commutativity for finite quotients $W^{J}$ to the $J^{c}$-dominant positions 
for which no game results in positive numbers on adjacent nodes.  
We then use a result of Stembridge to classify these ``adjacency-free'' 
positions. 
In \S \RootsNum\ we say precisely when all  
positive roots in the root system for a geometric 
representation of a finite Coxeter group can be obtained from a legal 
play sequence. 
In \S \DynkinClassNum\ we show that playing 
from a dominant position, the game will terminate if and only if it 
is played on a graph corresponding to a finite Coxeter group. 
(Another proof of this result based on ideas from \cite{ErikssonThesis} 
is given in \cite{DE}.) 
The geometric representations which connect the numbers game and 
Coxeter groups were introduced in \cite{ErikssonThesis} and 
\cite{ErikssonDiscrete} and studied further in \S 4.1-4.3 of \cite{BB} 
and \cite{DonCoxeter}.   
Definitions and results about these representations which are needed 
here are given in \S \PrelimNum.  There we also record 
several key results that are used throughout \S 
\ConvergenceNum-\DynkinClassNum: Eriksson's Strong Convergence Theorem, 
Eriksson's 
Comparison Theorem,  and 
Eriksson's Reduced Word Result.

\noindent 
{\bf Acknowledgments}\ \ 
We thank  
John Eveland for stimulating discussions during his work on 
an undergraduate research project \cite{Eveland} that led to the 
question addressed in \S \DynkinClassNum.  
We thank 
Norman Wildberger for sharing his perspective on the numbers game, 
including his observation about the appearance of 
``positive root functionals'' in numbers games on E-GCM graphs with 
integer amplitudes.  
We thank Bob Proctor for pointing us in the 
direction of Eriksson's work. 
We also thank Kimmo Eriksson for providing us with a copy of his 
thesis and for many helpful communications during the preparation of 
this paper.

%==================================================================
%\newpage 
\vspace{1ex} 

\noindent
{\Large \bf \PrelimNum.\ \ Definitions and preliminary results}

\vspace{1ex} 
Fix a positive integer $n$ and a totally ordered set $I_{n}$ with $n$ 
elements (usually $I_{n} := \{1<\ldots<n\}$).  An 
{\em E-generalized Cartan matrix}  
or {\em E-GCM}\myfootnote{Motivation for terminology: 
E-GCMs with integer entries are just  
generalized Cartan matrices, which 
are the starting point for the study of 
Kac--Moody 
algebras: beginning with a GCM, one can 
write down a list of the defining relations for a Kac--Moody 
algebra as well as its associated Weyl group (\cite{Kac}, \cite{Kumar}).  
Here we use the 
modifier ``E'' because of 
the relationship between these matrices and the combinatorics of 
Eriksson's E-games. Eriksson uses ``E'' for edge;  
he also allows for ``N-games'' 
where, in addition, nodes can be weighted.} is  
an $n \times n$ matrix $M = (M_{ij})_{i,j \in I_{n}}$  
with real entries satisfying the requirements that each 
main diagonal matrix entry is 2, that all other matrix entries are 
nonpositive, that if a matrix entry $M_{ij}$ is nonzero then its 
transpose entry $M_{ji}$ is also nonzero, and that if 
$M_{ij}M_{ji}$ is nonzero then $M_{ij}M_{ji} \geq 4$ or 
$M_{ij}M_{ji} = 4\cos^{2}(\pi/k_{ij})$ for some integer $k_{ij} \geq 
3$.  
These peculiar constraints on products of transpose pairs of matrix 
entries are precisely those required in order to guarantee ``strong 
convergence'' 
for E-games, 
cf.\ \StrongConvergenceTheorem\ below, Theorem 3.6 of 
\cite{ErikssonThesis}, Theorem 3.1 of 
\cite{ErikssonEur}. 
To an $n \times n$ E-generalized Cartan matrix 
$M = (M_{ij})_{i,j \in I_{n}}$ we associate a finite 
graph $\Gamma$ (which has undirected edges, 
no loops, and no multiple edges) 
as follows:     
The nodes $(\gamma_{i})_{i \in I_{n}}$ of $\Gamma$ are indexed 
by the set $I_{n}$, 
and   an edge is placed between nodes $\gamma_{i}$ and $\gamma_{j}$ 
if and only if $i \not= j$ 
and the matrix entries $M_{ij}$ and $M_{ji}$ are nonzero.  We call the pair 
$(\Gamma,M)$ an {\em E-GCM graph}. 
We depict a generic two-node E-GCM graph as follows:  

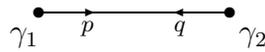
\vspace*{-0.075in}
\noindent
\begin{center}
\TwoCitiesGraphWithLabels
\end{center}

\vspace*{-0.075in}
\noindent
In this graph, $p = -M_{12}$  
and $q = -M_{21}$.  
We use \EGCMGraphCirclem 
for the collection of all two-node E-GCM 
graphs for which $M_{12}M_{21} = pq = 4\cos^{2}(\pi/m)$ 
for an integer $m > 3$; we use $m = \infty$ if $M_{12}M_{21} = pq \geq 
4$.  When $m = 3$ (i.e.\ $pq = 1$), 
we use an unlabelled edge \ATwoGraphNoEdgeLabels\ .  
An {\em E-Coxeter graph} will be any E-GCM graph whose 
connected components come from one of the collections of 
\ECoxeterGraphFigure. 

\begin{figure}[t]
\begin{center}
\ECoxeterGraphFigure: Families of connected E-Coxeter graphs. 
\vspace*{0.04in}

{\footnotesize (For adjacent nodes, the notation$\!\!$ 
\CircleIntegerm$\!\!$ means 
that the amplitude product on the edge is $4\cos^{2}(\pi/m)$;\\ for an 
unlabelled edge take $m=3$.  The asterisks for $\mathcal{E}_{6}$, 
$\mathcal{E}_{7}$, and $\mathcal{H}_{3}$ pertain to 
\AdjacencyFreeClassification.)}
\end{center}

\vspace*{-0.45in}
\begin{tabular}{cl}
$\mathcal{A}_{n}$ ($n \geq 1$) & \AnEGraph\\

$\mathcal{B}_{n}$ ($n \geq 3$) & \BnEGraph\\

$\mathcal{D}_{n}$ ($n \geq 4$) & \DnEGraph\\

$\mathcal{E}_{6}$ & \ESixEGraph\\

$\mathcal{E}_{7}$ & \ESevenEGraph\\

$\mathcal{E}_{8}$ & \EEightEGraph\\

$\mathcal{F}_{4}$ & \FFourEGraph\\

$\mathcal{H}_{3}$ & \HThreeEGraph\\

$\mathcal{H}_{4}$ & \HFourEGraph\\

$\mathcal{I}_{2}^{(m)}$ ($4 \leq m < \infty$) & \ITwoEGraph
\end{tabular}
\end{figure}

For the remainder of the paper the notation $(\Gamma,M)$ refers to an 
arbitrarily fixed E-GCM graph with nodes indexed by $I_{n}$,  
unless $(\Gamma,M)$ is otherwise specified.  
A {\em position} $\lambda = (\lambda_i)_{i \in 
I_{n}}$ is an assignment of real numbers to the nodes of 
$(\Gamma,M)$. 
The position $\lambda$ is 
{\em dominant} (respectively, {\em strongly dominant}) if 
$\lambda_{i} \geq 0$ 
(respectively $\lambda_i > 0$) for all $i \in I_{n}$; 
$\lambda$ is {\em nonzero} if at least one $\lambda_i \not= 0$. 
For $i \in I_{n}$, the {\em 
fundamental position} $\omega_i$ is the assignment of the number 
$1$ at node $\gamma_{i}$ and the number $0$ at all other nodes.  
Given a position $\lambda$ for $(\Gamma,M)$, to 
{\em fire} a node $\gamma_{i}$ is to change the number at each node 
$\gamma_{j}$ of $\Gamma$ by the transformation  
\[\lambda_j	 \longmapsto \lambda_j - 
M_{ij}\lambda_i,\] provided the number at node 
$\gamma_{i}$ is positive. Otherwise, node $\gamma_{i}$ is not allowed 
to be fired. 
In view of this transformation we think of entries of the E-GCM as 
{\em amplitudes}, and we sometimes refer to E-GCMs as 
{\em amplitude matrices}.  
The {\em numbers game} 
is the one-player game on $(\Gamma,M)$ in which the player 
(1) Assigns an initial position  
to the nodes of $\Gamma$; (2) Chooses a node with a positive 
number and fires the node to obtain a new position; and (3) 
Repeats step (2) for the new position if there is at least one node 
with a positive number.\myfootnote{Mozes 
studied numbers games on E-GCM graphs with integer amplitudes and for 
which the amplitude matrix $M$ is {\em symmetrizable} 
(i.e.\ there is a nonsingular  
diagonal matrix $D$ such that $D^{-1}M$ is symmetric). 
In \cite{Mozes} he obtained strong convergence results and a 
geometric characterization of the initial positions for which the 
game terminates.}  
Consider now the E-Coxeter graph in the $\mathcal{I}_{2}^{(4)}$ 
family depicted in \BtwoMarsDemo.  
As we can see in \BtwoMarsDemo, 
the numbers game terminates in a finite number of steps for any 
initial position and any legal sequence of node firings, 
if it is understood that the player  
will continue to fire as long as there is at least one 
node with a positive number.  In general, 
given a position $\lambda$, a {\em game sequence 
for} $\lambda$ is the (possibly empty, possibly infinite) sequence 
$(\gamma_{i_{1}}, \gamma_{i_{2}},\ldots)$, where 
$\gamma_{i_{j}}$ 
is the $j$th node that is fired in some 
numbers game with initial position $\lambda$.  
More generally, a {\em firing sequence} from some position $\lambda$ is an 
initial portion of some game sequence played from $\lambda$. The 
phrase {\em legal firing sequence} is used to emphasize that all node 
firings in the sequence are known or assumed to be possible. 
Note that a game sequence 
$(\gamma_{i_{1}}, 
\gamma_{i_{2}},\ldots,\gamma_{i_{l}})$ 
is of finite length $l$ 
(possibly with $l = 0$) if 
the number is nonpositive at each node after the $l$th firing. In 
this case we say the game sequence is {\em convergent} and the 
resulting position is the {\em terminal position} for the game 
sequence.  

\begin{figure}[t]
\begin{center}
\BtwoMarsDemo: The numbers game for an E-Coxeter graph in the 
$\mathcal{I}_{2}^{(4)}$ family. 

\vspace*{0.25in}
\setlength{\unitlength}{0.6in}
\begin{picture}(1.5,5.15) 
\put(0.25,5.1){\BTwoGraphForFigure}
\put(0.25,5.3){\scriptsize $a$}
\put(1.35,5.3){\scriptsize $b$}
\put(0.3,4.9){\vector(-4,-3){0.7}}
\put(1.4,4.9){\vector(4,-3){0.7}}
%==================
\put(1.75,3.8){\BTwoGraphForFigure}
\put(1.5,4.0){\scriptsize $a+2b$}
\put(2.75,4.0){\scriptsize $-b$}
\put(2.4,3.6){\vector(0,-1){0.6}}
%==================
\put(1.75,2.5){\BTwoGraphForFigure}
\put(1.45,2.7){\scriptsize $-a-2b$}
\put(2.75,2.7){\scriptsize $a+b$}
\put(2.4,2.3){\vector(0,-1){0.6}}
%==================
\put(1.75,1.2){\BTwoGraphForFigure}
\put(1.75,1.4){\scriptsize $a$}
\put(2.65,1.4){\scriptsize $-a-b$}
\put(2.1,1.0){\vector(-4,-3){0.7}}
%==================
\put(-1.25,3.8){\BTwoGraphForFigure}
\put(-1.35,4.0){\scriptsize $-a$}
\put(-0.25,4.0){\scriptsize $a+b$}
\put(-0.7,3.6){\vector(0,-1){0.6}}
%==================
\put(-1.25,2.5){\BTwoGraphForFigure}
\put(-1.5,2.7){\scriptsize $a+2b$}
\put(-0.35,2.7){\scriptsize $-a-b$}
\put(-0.7,2.3){\vector(0,-1){0.6}}
%==================
\put(-1.25,1.2){\BTwoGraphForFigure}
\put(-1.55,1.4){\scriptsize $-a-2b$}
\put(-0.15,1.4){\scriptsize $b$}
\put(-0.4,1.0){\vector(4,-3){0.7}}
%==================
\put(0.25,-0.1){\BTwoGraphForFigure}
\put(0.15,0.1){\scriptsize $-a$}
\put(1.25,0.1){\scriptsize $-b$}
%==================
\end{picture}
\end{center}

\vspace*{-0.25in}
\end{figure}

The preliminary results through \SubgraphLemma\ below also appear in 
\cite{DE} in the context of E-GCM graphs with integer 
amplitudes for use in a combinatorial 
proof of a result related to \EGCMTheorem.  
Proofs or references for 
these results are also given here.  
Following \cite{ErikssonThesis} and \cite{ErikssonEur}, we say 
the numbers game on an E-GCM graph $(\Gamma,M)$ is {\em strongly 
convergent} if given any initial position, every game sequence 
either diverges or converges to the same terminal position in the 
same number of steps. The next result  
follows from  
Theorem 3.1 of \cite{ErikssonEur} 
(or see Theorem 3.6 of \cite{ErikssonThesis}).  

\noindent
{\bf \StrongConvergenceTheorem\ (Eriksson's Strong Convergence Theorem)}\ \ 
{\sl The numbers game on a connected E-GCM graph 
is strongly 
convergent.}  

The 
following weaker result also applies when the 
E-GCM graph is not connected:  

\noindent
{\bf \StrongConvergenceCorollary} \ \ {\sl For any E-GCM graph, 
if a game sequence for an initial position  
$\lambda$ diverges, then all game  
sequences for $\lambda$ diverge.}
 
The next result is an immediate consequence of Theorem 4.3 of 
\cite{ErikssonThesis} or Theorem 4.5 of 
\cite{ErikssonDiscrete}.  
Eriksson's proof of this result in \cite{ErikssonThesis} uses only 
combinatorial and linear algebraic methods. 
 
\noindent 
{\bf \ComparisonTheorem\ (Eriksson's Comparison Theorem)}\ \ 
{\sl Given an E-GCM graph, suppose that a game sequence 
for an initial position $\lambda = (\lambda_{i})_{i \in 
I_{n}}$ converges.  Suppose that a position $\lambda' := (\lambda'_{i})_{i 
\in I_{n}}$ has the property that $\lambda'_{i} \leq 
\lambda_{i}$ for all $i \in I_{n}$.  Then some game sequence 
for the initial position $\lambda'$ also converges.}

Let $r$ be a positive real number.  
Observe that if 
$(\gamma_{i_{1}},\ldots,\gamma_{i_{l}})$ 
is a convergent 
game sequence for an initial position $\lambda = 
(\lambda_{i})_{i \in I_{n}}$, then 
$(\gamma_{i_{1}},\ldots,\gamma_{i_{l}})$ 
is a convergent 
game sequence for the initial position $r\lambda := 
(r\lambda_{i})_{i \in I_{n}}$. 
This observation and \ComparisonTheorem\ imply the following 
result: 

\noindent 
{\bf \ComparisonCorollary}\ \ {\sl 
Let $\lambda = (\lambda_{i})_{i \in 
I_{n}}$ be a dominant initial position such that $\lambda_{j} > 0$ 
for some $j \in I_{n}$. Suppose that a game sequence for $\lambda$ 
converges.  
Then some game   
sequence for the fundamental position $\omega_{j}$ also 
converges.}

The following is an immediate consequence of \ComparisonResults: 

\noindent 
{\bf \NotMarsFriendlyLemma}\ \ {\sl An E-GCM graph is not admissible if 
for each fundamental position there is a divergent game  
sequence.} 

The following is 
proved easily with an induction argument on the number of 
nodes.  

\noindent 
{\bf \EveryNodeFiredLemma}\ \ {\sl Suppose $(\Gamma,M)$ is connected 
with nonzero dominant position $\lambda$. 
Then in any 
convergent 
game sequence for $\lambda$, every node of $\Gamma$ is fired at least 
once.}

If $I'_{m}$ is a subset of the node set $I_{n}$ of 
$(\Gamma,M)$, 
then let $\Gamma'$ be the subgraph of $\Gamma$ with node set $I'_{m}$ and 
the induced set of edges, and let $M'$ be the corresponding 
submatrix 
of the amplitude matrix $M$. We call $(\Gamma',M')$ an {\em E-GCM 
subgraph} of $(\Gamma,M)$.  In light of \NewLemmaList, the 
following result amounts to an observation. 

\noindent 
{\bf \SubgraphLemma}\ \ {\sl If a connected E-GCM graph is 
admissible, then any 
connected E-GCM subgraph is also admissible.} 

Define the associated Coxeter group 
$W = W(\Gamma,M)$ to be the Coxeter group with identity 
denoted $\varepsilon$, 
generators $\{s_{i}\}_{i \in I_{n}}$, and defining relations $s_{i}^{2} = 
\varepsilon$ for $i \in I_{n}$ and 
$(s_{i}s_{j})^{m_{ij}} = \varepsilon$ for all $i \not= j$, where the 
$m_{ij}$ are determined as follows: 

\vspace*{-0.1in}
\[m_{ij} = \left\{\begin{array}{cl}
k_{ij} & \mbox{\hspace*{0.25in} if 
$M_{ij}M_{ji} = 4\cos^{2}(\pi/k_{ij})$ for some integer $k_{ij} \geq 2$}\\ 
\infty & \mbox{\hspace*{0.25in} if 
$M_{ij}M_{ji} \geq 4$} 
\end{array}\right.\]  
(Conventionally, $m_{ij} = \infty$ means there is no relation between 
generators 
$s_{i}$ and $s_{j}$.) 
Throughout the paper, $W$ denotes the Coxeter group $W(\Gamma,M)$ 
associated to an arbitrarily fixed E-GCM graph $(\Gamma,M)$ with index 
set $I_{n}$. 
One can think of the E-GCM graph as a refinement of the information 
from the Coxeter graph for the associated Coxeter group.    
Observe that any Coxeter group on a finite set of generators is isomorphic 
to the Coxeter group associated to some E-GCM graph. 
The Coxeter group $W$ is {\em irreducible} if $\Gamma$ is connected.   
Let ${\myl}$ denote the length function for 
the $W$. 
An expression $s_{i_{p}}{\cdots}s_{i_{2}}s_{i_{1}}$ for an element of 
$W$ is {\em reduced} 
if $\myl(s_{i_{p}}{\cdots}s_{i_{2}}s_{i_{1}}) = p$. An empty product 
in $W$ is taken as $\varepsilon$. 
For a firing sequence $(\gamma_{i_{1}}, 
\gamma_{i_{2}}, \ldots, \gamma_{i_{p}})$ from some initial position 
on $(\Gamma,M)$, 
the corresponding element of $W$ is taken to be  
$s_{i_{p}}\cdots{s}_{i_{2}}s_{i_{1}}$. 
Parts (1) and (2) of what we call Eriksson's Reduced Word Result 
follow respectively from Propositions 4.1 and 4.2 of 
\cite{ErikssonDiscrete}. 

\noindent 
{\bf \ErikssonWordProposition\ (Eriksson's 
Reduced Word Result)}\ \ {\sl (1) If 
$(\gamma_{i_{1}}, 
\gamma_{i_{2}}, \ldots, \gamma_{i_{p}})$ is a legal sequence of node 
firings in a numbers game played  
from some initial position on $(\Gamma,M)$, 
then $s_{i_{p}}\cdots{s}_{i_{2}}s_{i_{1}}$ is a reduced expression for 
the corresponding element of $W$. (2) If 
$s_{i_{p}}\cdots{s}_{i_{2}}s_{i_{1}}$ is a reduced expression for an 
element of $W$, then $(\gamma_{i_{1}}, 
\gamma_{i_{2}}, \ldots, \gamma_{i_{p}})$ is a legal sequence of node 
firings in a numbers game played  
from any given strongly dominant position on $(\Gamma,M)$.} 

To conclude this section 
we summarize results from \cite{DonCoxeter} concerning certain 
geometric representations of Coxeter groups introduced by Eriksson in 
\cite{ErikssonThesis} and \cite{ErikssonDiscrete}.  
Let $V$ be a real $n$-dimensional vector space freely generated by 
$(\alpha_{i})_{i \in I_{n}}$ (elements of this ordered basis are {\em 
simple roots}).  
Equip $V$ with a possibly asymmetric bilinear form $B: V \times V 
\rightarrow \mathbb{R}$ defined on the 
basis $(\alpha_{i})_{i \in I_{n}}$ by 
$B(\alpha_{i},\alpha_{j}) := \frac{1}{2}M_{ij}$. 
For each $i \in I_{n}$ define an operator 
$S_{i}: V \rightarrow V$ by the rule $S_{i}(v) := v - 
2B(\alpha_{i},v)\alpha_{i}$ for each $v \in V$.  One can check 
that $S_{i}^{2}$ is the identity 
transformation, so $S_{i} \in GL(V)$.  

As can be seen for example in 
\cite{BB} Theorem 4.2.2, there is a unique homomorphism 
$\sigma_{M}: W \rightarrow GL(V)$ for which $\sigma_{M}(s_{i}) = 
S_{i}$.  
Theorem 4.2.7 of \cite{BB} shows that 
$\sigma_{M}$ is injective. 
We call $\sigma_{M}$ a {\em geometric representation} 
of $W$. 
We now have $W$ acting on $V$, and for all $w \in W$ and $v \in V$ we 
write $w.v$ for $\sigma_{M}(w)(v)$. 
Define $\Phi_{M} := \{\alpha \in V\, |\, \alpha = w.\alpha_{i} \mbox{ for 
some } i \in I_{n} \mbox{ and } w \in W\}$.  For each $w \in W$,  
$\sigma_{M}(w)$ permutes $\Phi_{M}$, so $\sigma_{M}$ induces an 
action of $W$ on $\Phi_{M}$. Evidently,  
$\Phi_{M} = -\Phi_{M}$.  Elements of $\Phi_{M}$ are 
{\em roots} and are necessarily nonzero.  
If $\alpha = \sum c_{i}\alpha_{i}$ is a root with all 
$c_{i}$ nonnegative (respectively nonpositive), then say $\alpha$ is a {\em 
positive} (respectively {\em negative}) root.  
Let $\Phi_{M}^{+}$ and $\Phi_{M}^{-}$ 
denote the collections of positive and negative roots respectively.  
Let $w \in W$ and $i \in I_{n}$.  
Proposition 4.2.5 of \cite{BB} states: {\sl If 
$\myl(ws_{i}) > 
\myl(w)$ then $w.\alpha_{i} \in \Phi_{A}^{+}$, and  
if $\myl(ws_{i}) < 
\myl(w)$ then $w.\alpha_{i} \in \Phi_{A}^{-}$.} 
It follows that $\Phi_{M}$ is 
partitioned by $\Phi_{M}^{+}$ and $\Phi_{M}^{-}$. 

We say two adjacent nodes $\gamma_{i}$ and $\gamma_{j}$ 
in $(\Gamma,M)$ are {\em odd-neighborly} if $m_{ij}$ is odd, {\em 
even-neighborly} if $m_{ij} \geq 4$ is even, and $\infty$-{\em 
neighborly} if 
$m_{ij} = \infty$.  
When 
$m_{ij}$ is odd and 
$M_{ij} \not= M_{ji}$, we say that the adjacent nodes $\gamma_{i}$ 
and $\gamma_{j}$ form an {\em odd asymmetry}. 
For odd $m_{ij}$, 
let $v_{ji}$ be the element 
$(s_{i}s_{j})^{(m_{ij}-1)/2}$, and 
set $K_{ji} := \frac{-M_{ji}}{2\cos(\pi/m_{ij})}$, which is positive. 
It is a consequence of \CoxeterTwoGeneratorAnalysis\ of 
\cite{DonCoxeter} that  
$v_{ji}.\alpha_{i} = K_{ji}\alpha_{j}$.  
Observe that $K_{ij}K_{ji} = 1$ and moreover that 
$v_{ij} = v_{ji}^{-1}$.  
A {\em path of odd neighbors} (or {\em ON-path}, 
for short) 
in $(\Gamma,M)$ is a sequence  
$\mathcal{P} := 
[\gamma_{i_{0}},\gamma_{i_{1}},\ldots,\gamma_{i_{p}}]$ of 
nodes from $\Gamma$ for which consecutive pairs are odd-neighborly.  
This ON-path has length $p$, and we allow 
ON-paths to have length zero.  
We say $\gamma_{i_{0}}$ and 
$\gamma_{i_{p}}$ are the {\em start} and {\em end} nodes of the ON-path, 
respectively.  
Let $w_{_{\mathcal{P}}} \in W$ be the Coxeter group element 
$v_{i_{p}i_{p-1}}\cdots{v}_{i_{2}i_{1}}v_{i_{1}i_{0}}$, and 
let $\Pi_{_{\mathcal{P}}} := 
K_{i_{p}i_{p-1}}\cdots{K}_{i_{2}i_{1}}K_{i_{1}i_{0}}$, where 
$w_{_{\mathcal{P}}} = \varepsilon$ with $\Pi_{_{\mathcal{P}}} = 1$ 
when $\mathcal{P}$ has length zero. 
Note that $w_{_{\mathcal{P}}}.\alpha_{i_{0}} = 
\Pi_{_{\mathcal{P}}}\alpha_{i_{p}}$. The next result follows from  
\CoxeterHumphreysTheorem\ of \cite{DonCoxeter}. 

\noindent 
{\bf \HumphreysTheorem}\ \ 
{\sl Let $w \in W$ and $i \in I_{n}$. 
(1) 
Then $w.\alpha_{i} = 
K\alpha_{x}$ for some $x \in I_{n}$ and some $K > 0$  
if and only if $w.\alpha_{i} = 
w_{_{\mathcal{P}}}.\alpha_{i}$ for some ON-path $\mathcal{P} = 
[\gamma_{i_{0}=i},
\gamma_{i_{1}},\ldots,\gamma_{i_{p-1}},\gamma_{i_{p}=x}]$, in which 
case 
$K = \Pi_{_{\mathcal{P}}}$.  
(2) 
Similarly $w.\alpha_{i} = 
K\alpha_{x}$ for some $x \in I_{n}$ and some $K<0$  
if and only if $w.\alpha_{i} = (w_{_{\mathcal{P}}}s_{i}).\alpha_{i}$ 
for some ON-path $\mathcal{P} = 
[\gamma_{i_{0}=i},
\gamma_{i_{1}},\ldots,\gamma_{i_{p-1}},\gamma_{i_{p}=x}]$, in which 
case 
$K = -\Pi_{_{\mathcal{P}}}$.} 

An ON-path 
$\mathcal{P} = [\gamma_{i_{0}},\ldots,\gamma_{i_{p}}]$ 
is an {\em ON-cycle} if $\gamma_{i_{p}} = \gamma_{i_{0}}$.  
It is a {\em unital} ON-cycle if 
$\Pi_{\mathcal{P}} = 1$.  For ON-paths $\mathcal{P}$ and $\mathcal{Q}$, 
write $\mathcal{P} \sim 
\mathcal{Q}$ and say $\mathcal{P}$ and $\mathcal{Q}$ are {\em 
equivalent}  
if these ON-paths have 
the same start and end nodes and 
$\Pi_{_{\mathcal{P}}} = \Pi_{_{\mathcal{Q}}}$. This is an 
equivalence relation on the set of all ON-paths. 
An ON-path $\mathcal{P}$ 
is {\em simple} if it has no repeated nodes with the possible 
exception that the start and end nodes may coincide.  
We say  
$(\Gamma,M)$ is {\em unital ON-cyclic} if and only if 
$\Pi_{_{\mathcal{C}}} = 1$ for all ON-cycles $\mathcal{C}$.  
Note that $(\Gamma,M)$ is unital ON-cyclic if and only if 
$\mathcal{P} \sim \mathcal{Q}$ whenever 
$\mathcal{P}$ and $\mathcal{Q}$ are ON-paths 
with the same start and end nodes.  
The property that $(\Gamma,M)$ has no odd asymmetries is sufficient 
but not necessary to imply that $(\Gamma,M)$ is unital ON-cyclic. 
An E-GCM graph is {\em ON-connected} if 
any two nodes can be joined by an ON-path.  An {\em ON-connected 
component} of $(\Gamma,M)$ is an E-GCM subgraph 
$(\Gamma',M')$ whose nodes 
form a maximal collection of nodes in $(\Gamma,M)$ which 
can be pairwise joined by ON-paths.  
For any $\alpha \in \Phi_{M}$, 
set $\mathfrak{S}_{M}(\alpha) := 
\{K\alpha|K\in\mathbb{R}\}\cap\Phi_{M}^{+}$. 
The next result is \CoxeterTFAE\ 
of \cite{DonCoxeter}. 

\noindent 
{\bf \TFAE}\ \ 
{\sl 
Choose any ON-connected component $(\Gamma',M')$ of 
$(\Gamma,M)$, and let $J := \{x \in I_{n}\}_{\gamma_{x} \in \Gamma'}$. 
Then $(\Gamma',M')$ is unital ON-cyclic if and only if 
$|\mathfrak{S}_{M}(\alpha_{x})| < 
\infty$ for some $x \in J$ if and only if $|\mathfrak{S}_{M}(\alpha_{x})| < 
\infty$ for all $x \in J$, in which case 
we have $|\mathfrak{S}_{M}(\alpha_{x})| = 
|\mathfrak{S}_{M}(\alpha_{y})|$ for all $x, y \in 
J$.} 

%===================
% Formula for $N_{M}(w)$
%=================== 
For any $w \in W$, set $N_{M}(w) := 
\{\alpha \in \Phi_{M}^{+}\, |\, w.\alpha \in \Phi_{M}^{-}\}$. The 
following is \CoxeterPositiveToNegativeRootsLemma\ of \cite{DonCoxeter}. 

\noindent
{\bf \PositiveToNegativeRootsLemma}\ \ {\sl For any $i \in I_{n}$, 
$s_{i}(\Phi_{M}^{+}\setminus\mathfrak{S}_{M}(\alpha_{i})) = 
\Phi_{M}^{+}\setminus\mathfrak{S}_{M}(\alpha_{i})$.  Now let $w \in 
W$.  
If $w.\alpha_{i} \in \Phi_{M}^{+}$, then $N_{M}(ws_{i}) = s_{i}(N_{M}(w)) 
\disjointunion 
\mathfrak{S}_{M}(\alpha_{i})$, a disjoint union. 
If $w.\alpha_{i} \in \Phi_{M}^{-}$, then $N_{M}(ws_{i}) = 
s_{i}(N_{M}(w)\setminus\mathfrak{S}_{M}(\alpha_{i}))$.} 

When $(\Gamma,M)$ is ON-connected and unital ON-cyclic, let 
$f_{\Gamma,M} := |\mathfrak{S}_{M}(\alpha_{x})|$ for 
any fixed $x \in I_{n}$.  
For $J \subseteq I_{n}$, 
let $\mathfrak{C}(J)$ denote the set of all 
ON-connected components of $(\Gamma,M)$ containing  
some node from the set $\{\gamma_{x}\}_{x \in J}$. 
The next result is 
\CoxeterLengthProposition\ of \cite{DonCoxeter}. 

\noindent 
{\bf \LengthProposition}\ \ 
{\sl Let $w \in W$ with $p = \myl(w) > 0$.  
(1) Then $N_{M}(w)$ is 
finite if and only if $w$ has a reduced expression 
$s_{i_{1}}{\cdots}s_{i_{p}}$ for which 
$\mathfrak{S}_{M}(\alpha_{i_{q}})$ is finite for all $1 \leq q \leq p$ 
if and only if every 
reduced expression $s_{i_{1}}{\cdots}s_{i_{p}}$ for $w$ has 
$\mathfrak{S}_{M}(\alpha_{i_{q}})$ finite for all $1 \leq q \leq p$. 
(2) Now suppose $w = s_{i_{1}}{\cdots}s_{i_{p}}$ and 
$N_{M}(w)$ is finite.  Let $J := \{i_{1},\ldots,i_{p}\}$. 
In view of (1), let 
$f_{1}$ be the min and $f_{2}$ the max of all integers in the set 
$\{f_{\Gamma',M'}\, |\, (\Gamma',M') \in \mathfrak{C}(J)\}$.  
Then $f_{1}\, \myl(w) \leq |N_{M}(w)| \leq f_{2}\, \myl(w)$.} 
 
We have the natural pairing $\langle \lambda, v \rangle := 
\lambda(v)$ for elements $\lambda$ in the dual space $V^{*}$ and 
vectors $v$ in $V$. 
We think of $V^{*}$ as the space of positions for 
numbers games played on $(\Gamma,M)$:   
For $\lambda \in V^{*}$, the numbers for the corresponding 
position are $(\lambda_{i})_{i \in I_{n}}$ where for each $i \in 
I_{n}$ we have $\lambda_{i} := \langle \lambda, \alpha_{i} 
\rangle$.  
Regard the fundamental positions $(\omega_{i})_{i \in I_{n}}$ to be the 
basis for $V^{*}$ dual to the basis $(\alpha_{j})_{j \in I_{n}}$ for 
$V$ relative to the natural pairing $\langle \cdot,\cdot \rangle$, so $\langle 
\omega_{i}, \alpha_{j} \rangle = \delta_{ij}$. 
Given $\sigma_{M}: W \rightarrow GL(V)$, the 
contragredient representation 
$\sigma_{M}^{*}: W \rightarrow GL(V^{*})$ 
is determined by $\langle \sigma_{M}^{*}(w)(\lambda), v \rangle = 
\langle \lambda, \sigma_{M}(w^{-1})(v) \rangle$.  From here on, when $w 
\in W$ and $\lambda \in V^{*}$, 
write $w.\lambda$ for $\sigma_{M}^{*}(w)(\lambda)$.  
Then $s_{i}.\lambda$ is the result of firing node $\gamma_{i}$ when 
the E-GCM graph is assigned position $\lambda$, whether the firing is 
legal or not. 
We have a one-to-one correspondence between roots and 
certain elements of $V^{**}$:
Given a root $\alpha$, the {\em root functional} 
$\phi_{\alpha}: V^{*} \rightarrow \mathbb{R}$ is 
given by $\phi_{\alpha}(\mu) = \langle \mu, \alpha \rangle$, and 
$\phi_{\alpha}$ is {\em positive} (resp.\ {\em negative}) if $\alpha 
\in \Phi_{M}^{+}$ (resp.\ $\Phi_{M}^{-}$). 

\noindent 
{\bf \EquivalenceRemark}\ \ 
From the definitions one sees that the following 
are equivalent: 
(1) $(\gamma_{i_{1}},\ldots,\gamma_{i_{p}})$ is legally played 
from some position $\lambda$, 
(2) $\langle s_{i_{q-1}}\cdots{s_{i_{1}}}.\lambda, \alpha_{i_{q}} 
\rangle > 0$ for $1 \leq q \leq p$,  
(3) $\langle \lambda, \beta_{q} 
\rangle > 0$ where $\beta_{q} := 
s_{i_{1}}\cdots{s_{i_{q-1}}}.\alpha_{i_{q}}$ for $1 \leq q \leq p$, 
(4) $\phi_{\beta_{q}}(\lambda) > 0$ for $1 \leq q \leq p$. 
That 
$\beta_{q} \in \Phi_{M}^{+}$ for $1 \leq q \leq p$ follows from 
\cite{BB} Proposition 4.2.5 and the fact that 
$\myl(s_{i_{1}}\cdots{s}_{i_{q-1}}) < 
\myl(s_{i_{1}}\cdots{s}_{i_{q-1}}s_{i_{q}})$.\hfill\QED

Let $D$ be the set of dominant positions. 
The {\em Tits 
cone} is $U_{M} := \cup_{w \in W}wD$.  The next result is 
\CoxeterHRTResult\ of \cite{DonCoxeter}. 

\noindent 
{\bf \HRTResult}\ \ {\sl Suppose $(\Gamma,M)$ is connected and unital 
ON-cyclic.  If the Coxeter group $W$ 
is infinite, then 
$U_{M} \cap (-U_{M}) = \{0\}$.} 
 
In \S 4 of \cite{ErikssonDiscrete}, Eriksson  
characterizes the set of initial positions for which the game 
converges. In contrast to \cite{ErikssonDiscrete}, 
here we fire at nodes with positive rather than negative 
numbers, so we have $-U_{M}$ instead of $U_{M}$ in the following 
statement.  

\noindent 
{\bf \TitsConeConvergenceResult\ (Eriksson)}\ \ 
{\sl The set of initial positions for which the numbers game on the 
E-GCM graph $(\Gamma,M)$ converges is precisely $-U_{M}$.}

%==================================================================
\vspace{1ex} 

\noindent
{\Large \bf \ConvergenceNum.\ \ 
Extensions of Eriksson's Reduced Word Result for dominant 
positions}

\vspace{1ex} 
In this section we consider legal play sequences from dominant 
positions with a specified set $J$ of nodes where the numbers are zero.  
This leads to 
certain extensions of Eriksson's Reduced Word Result 
in \KeyWJResult\  and \KeyWJResultCorollary.  
Eriksson's Strong Convergence Theorem 
is used in deriving two 
corollaries to \KeyWJResult.   
For any $J \subseteq I_{n}$, 
$W_{J}$ is the subgroup generated by 
$\{s_{i}\}_{i \in J}$, a {\em parabolic} subgroup, and $W^{J} := \{w 
\in W\, |\, \myl(ws_{j}) > \myl(w) \mbox{ for all } j \in J\}$  
is the set of {\em minimal coset representatives} (see \cite{BB} Ch.\ 
2). 
When $J = \emptyset$, $W_{J}$ is the one-element group and $W^{J} = W$.  
If $W$ is finite,  
we may choose the (unique) longest 
element $w_{0}$ in $W$.  
Since we must have $\myl(w_{0}s_{i}) < 
\myl(w_{0})$ for all $i \in I_{n}$, it follows that $w_{0}.\alpha_{i} 
\in \Phi_{M}^{-}$ for all $i$.  So if $\alpha = \sum c_{i}\alpha_{i} \in \Phi_{M}^{+}$, 
then $w_{0}.\alpha \in \Phi_{M}^{-}$, i.e.\   
$N_{M}(w_{0}) = \Phi_{M}^{+}$.  
More generally, 
for any $W$ (not necessarily finite) and for any subset 
$J$ of $I_{n}$, we let $(w_{0})_{_{J}}$ denote 
the longest element of $W_{J}$ when  
$W_{J}$ is finite. 

%================
% A lemma needed to generalize Eriksson 
%================
\noindent 
{\bf \KeyWJResultLemma}\ \ 
{\sl   
Let $J \subseteq I_{n}$, 
and suppose 
$W_{J}$ is finite.  
Suppose $\alpha = \sum_{j \in 
J}c_{j}\alpha_{j}$ is a root in $\Phi_{M}^{+}$.  
Then 
$(w_{0})_{_{J}}.\alpha \in \Phi_{M}^{-}$.} 

{\em Proof.} 
Note that any element of $W_{J}$ preserves 
the subspace $V_{J} := \mathrm{span}_{\mathbb{R}}\{\alpha_{j}\}_{j \in 
J}$.  As seen just above, $(w_{0})_{_{J}}$ will send each simple 
root $\alpha_{j}$ for $j \in J$ to some root in $\Phi_{M}^{-}$.  
Then $(w_{0})_{_{J}}.\alpha \in \Phi_{M}^{-}$. \hfill\QED

%======================
% The $W^{J}$ result
%======================
In what follows, for any subset $J$ of 
$I_{n}$, a position $\lambda$ is $J^{c}$-{\em dominant} if its zeros 
are precisely on the nodes in set $J$, i.e.\  
$\lambda = \sum_{i \in 
I_{n}\setminus{J}}\lambda_{i}\omega_{i}$ with $\lambda_{i} > 0$ for 
all $i \in I_{n}\setminus{J}$. 
Part (2) of Eriksson's Reduced 
Word Result 
and the ``if'' direction of 
Theorem 4.3.1.{\em iv} of \cite{BB} are the $J = \emptyset$ case 
of our next result. 

\noindent 
{\bf \KeyWJResult}\ \ 
{\sl Let $J \subseteq I_{n}$ and let $\lambda$ be $J^{c}$-dominant. 
Suppose $W_{J}$ is finite. Let $s_{i_{p}}\cdots{s}_{i_{2}}s_{i_{1}}$ be any reduced 
expression for an element  
of $W^{J}$.  
Then 
$(\gamma_{i_{1}},\ldots,\gamma_{i_{p}})$ is a legal sequence of node 
firings from initial position $\lambda$.  That is, the root 
$s_{i_{1}}s_{i_{2}}\cdots{s_{i_{q-1}}}.\alpha_{i_{q}}$ is positive for $1 \leq 
q \leq p$.} 

{\em Proof.}  
By \EquivalenceRemark, 
we must show that 
$\langle \lambda ,\beta_{q} 
\rangle > 0$
for $1 \leq q \leq p$, where $\beta_{q} := 
s_{i_{1}}\cdots{s}_{i_{q-1}}.\alpha_{i_{q}}$.  
Suppose 
$s_{j_{r}}\cdots{s}_{j_{2}}s_{j_{1}}$ is a reduced expression for some 
$v_{_{J}} \in W_{J}$.  Since 
$s_{i_{p}}\cdots{s}_{i_{2}}s_{i_{1}}s_{j_{1}}\cdots{s}_{j_{r}}$
is reduced (cf.\ Proposition 2.4.4 of \cite{BB}), it follows that 
$\myl(v_{_{J}}s_{i_{1}}\cdots{s}_{i_{q-2}}s_{i_{q-1}}) <  
\myl(v_{_{J}}s_{i_{1}}\cdots{s}_{i_{q-1}}s_{i_{q}})$. 
In particular 
$v_{_{J}}.\beta_{q} \in \Phi_{M}^{+}$ for all $v_{_{J}} 
\in W_{J}$.  
We wish to show that 
$\beta_{q}$ cannot be 
contained in $\mathrm{span}_{\mathbb{R}}\{\alpha_{j}\}_{j \in J}$.  
Suppose otherwise, so  
$\beta_{q} = 
\sum_{j \in J}c_{j}\alpha_{j}$.  
\EquivalenceRemark\ shows that $\beta_{q} \in \Phi_{M}^{+}$ for $1 \leq q 
\leq p$.  
But now the finiteness 
of $W_{J}$ and \KeyWJResultLemma\ imply that  
$(w_{0})_{_{J}}.\beta_{q} \in \Phi_{M}^{-}$, a contradiction.  
Then it must be the case that 
$\beta_{q} = 
\sum_{i \in I_{n}}c_{i}\alpha_{i}$ with $c_{k} > 0$ for some $k \in 
I_{n}\setminus{J}$.  So 
$\langle 
\lambda, \beta_{q} \rangle 
= \langle \lambda , \sum_{i \in I_{n}}c_{i}\alpha_{i} \rangle 
= \sum_{i \in I_{n}}c_{i}\lambda_{i}$, which is positive since all 
$c_{i}$'s are nonnegative, $\lambda_{k} > 0$, and $c_{k} > 
0$.\hfill\QED  

It is an open question whether the finiteness hypothesis of 
\KeyWJResult\ can be relaxed. 
See \S \DynkinClassNum\ for comments on a possible connection between 
\KeyWJResult\ and \EGCMTheorem. 
Let $\mathfrak{P}(\lambda)$ denote the set of all positions obtainable 
from legal firing sequences in numbers games with initial position 
$\lambda$.  Clearly $\mathfrak{P}(\lambda) \subseteq W\lambda$, where 
the latter is the orbit of $\lambda$ under the $W$-action on 
$V^{*}$.  
Since the 
statement of  
Theorem 5.13 of \cite{HumCoxeter} holds for geometric 
representations, then $W_{J}$ is the full stabilizer of any 
$J^{c}$-dominant $\lambda$, 
so $W\lambda$ and $W^{J}$ can be identified. 
So from \KeyWJResult\ we see that for $J^{c}$-dominant $\lambda$ with  
$W_{J}$ finite, then $\mathfrak{P}(\lambda) = W\lambda$.  
The $J = \emptyset$ version of the previous statement is part ({\em 
ii}) of Theorem 4.3.1 of \cite{BB}. 

For finite $W$, we use $(w_{0})^{^{J}}$ to denote 
the minimal coset representative for $w_{0}W_{J}$. 

\noindent 
{\bf \LengthUsingLongestWord}\ \ 
{\sl Suppose $W$ is finite.  Let $J \subseteq I_{n}$.  
Then all game sequences for any $J^{c}$-dominant $\lambda$ have length 
$\myl((w_{0})^{^{J}}) = \myl(w_{0}) - \myl((w_{0})_{_{J}})$.} 

{\em Proof.} 
\KeyWJResult\ implies that there is a game sequence for $\lambda$ with 
length $\myl((w_{0})^{^{J}}) = \myl(w_{0}) - \myl((w_{0})_{_{J}})$.  
By Eriksson's 
Strong Convergence Theorem, this must be the length of any game 
sequence for $\lambda$.\hfill\QED

For finite Coxeter groups, the next result strengthens Part (1) of Eriksson's 
Reduced Word Result. 
At this time it is an open question whether the finiteness hypothesis 
for $W$ can be relaxed. 

\noindent 
{\bf \KeyWJResultCorollary}\ \ 
{\sl Let $J \subseteq I_{n}$ and let $\lambda$ be any $J^{c}$-dominant 
position. 
Suppose $W$ is finite. Suppose $\mathbf{s} := 
(\gamma_{i_{1}},\ldots,\gamma_{i_{p}})$ 
is a legal firing sequence for played from $\lambda$.  
Then $w := s_{i_{p}}\cdots{s}_{i_{2}}s_{i_{1}}$ is a reduced 
expression for an element  
of $W^{J}$.  Moreover, $\mathbf{s}$ is a game sequence if and only if 
$w = (w_{0})^{^{J}}$.}

{\em Proof.}  
By \LengthUsingLongestWord, we may 
extend the legal firing sequence $\mathbf{s}$ 
to some game sequence 
$\selt' := 
(\gamma_{i_{1}},\ldots,\gamma_{i_{p}},\gamma_{i_{p+1}},\ldots,\gamma_{i_{L}})$ 
with $L = \myl(w_{0}) - \myl((w_{0})_{_{J}}) \geq p$.  
Let 
$v := s_{i_{L}}\cdots{s}_{i_{p+2}}s_{i_{p+1}}$, 
and $u := vw$.  By Part (1) of Eriksson's Reduced Word 
Result, 
$w$, $v$, and $u$ are reduced.  In particular, 
$\myl(u) = L$.  Write $u = u^{^{J}}u_{_{J}}$ for $u^{^{J}} \in W^{J}$ and 
$u_{_{J}} \in W_{J}$.  By \KeyWJResult, we may take 
a legal firing sequence $\telt := (\gamma_{j_{1}},\ldots,\gamma_{j_{K}})$  
from $\lambda$ corresponding to some 
reduced expression for $u^{^{J}}$.  Now $u.\lambda$ is the terminal 
position for the game sequence $\selt'$ played from $\lambda$.  
Since $u.\lambda = u^{^{J}}.\lambda$, then $\telt$ is a game sequence 
terminating at this same position.  By Eriksson's Strong Convergence 
Theorem, it must be the case that 
$\myl(u^{^{J}}) = K = L = \myl(u)$.  Hence $u_{_{J}} = \varepsilon$ 
and $u = u^{^{J}} \in W^{J}$.  Now write $w = w^{^{J}}w_{_{J}}$ for 
$w^{^{J}} \in W^{J}$ and $w_{_{J}} \in W_{J}$.  If 
$w_{_{J}} \not= \varepsilon$, then $w_{_{J}}$ has a reduced 
expression ending in $s_{j}$ for some $j \in J$.  Then 
$\myl(us_{j}) = \myl(vws_{j}) = \myl(vw^{^{J}}w_{_{J}}s_{j}) < 
\myl(u)$.  But this contradicts the fact that $u \in W^{J}$.  Hence 
$w_{_{J}} = \varepsilon$, so $w = w^{^{J}} \in W^{J}$. 
By Proposition 2.4.4 of \cite{BB}, $\myl(u(w_{0})_{_{J}}) = L+\myl(\, 
(w_{0})_{_{J}})$.  Since $L+\myl(\, 
(w_{0})_{_{J}}) = \myl(w_{0})$, then $u(w_{0})_{_{J}} = w_{0} = 
(w_{0})^{^{J}}(w_{0})_{_{J}}$, so $u = (w_{0})^{^{J}}$. 
It now follows that $w = (w_{0})^{^{J}}$if and only if $\mathbf{s}$ is a game 
sequence.\hfill\QED

%==================================================================
\newpage

\noindent
{\Large \bf \AdjacencyFreeNum.\ \ 
Adjacency-free positions and full commutativity 
of Coxeter group 
elements}

\vspace{1ex} 
In this section we study dominant positions whose numbers games are 
all equivalent up to a notion of interchanging moves.  We say these 
positions are ``adjacency-free.''  For finite $W$, we classify the 
adjacency-free positions by showing how they correspond with 
quotients $W^{J}$ whose elements are fully commutative in the sense 
of \cite{StemFC} (see also 
\cite{StemFCEnumerate}; see \cite{FanJAC} and \cite{FanJAMS} 
for full commutativity in a different 
context).  Adjacency-free positions have other connections to the 
literature.  
In what follows, a Weyl group is a Coxeter group for which each 
$m_{ij} \in \{2,3,4,6,\infty\}$.  
In Proposition 3.1 of \cite{PrEur}, 
Proctor shows that for finite irreducible Weyl groups $W$, those 
quotients for which the Bruhat order $(W^{J},\leq)$ (see \cite{BB}) 
is a lattice have $|J^{c}| = 
1$ and correspond  
precisely to the adjacency-free fundamental positions for the 
connected ``Dynkin 
diagrams of finite type'' (E-Coxeter graphs with integer amplitudes).  
In Proposition 
3.2 of that paper, he shows 
that these lattices 
are, in fact, distributive. 
In \cite{DW} we use information obtained from numbers games played 
from adjacency-free fundamental positions  on Dynkin diagrams of 
finite type to construct certain ``fundamental'' posets.  We show that 
the distributive lattices of order ideals obtained from certain 
combinations of our fundamental posets can be used to produce Weyl 
characters and in some cases explicit constructions irreducible 
representations of the corresponding semisimple Lie algebra.  For 
rank two versions of these posets and distributive lattices, see 
\cite{ADLMPPW} and \cite{ADLP}. 
When an 
adjacency-free fundamental position for a Dynkin diagram of finite 
type corresponds 
to a ``minuscule'' fundamental weight (see \cite{PrEur}, 
\cite{PrDComplete}, 
\cite{StemFC}), then our fundamental poset coincides with 
the corresponding ``wave'' poset 
of \cite{PrDComplete} and  ``heap'' of 
\cite{StemFC}.

For a firing sequence $(\gamma_{i_{1}}, \gamma_{i_{2}}, \ldots)$ from 
a position $\lambda$, any position 
$s_{i_{j}}\cdots{s}_{i_{1}}.\lambda$ (including $\lambda$ itself) 
is an {\em intermediate 
position} for the sequence.  
A game sequence played from $\lambda$ is {\em adjacency-free} if no 
intermediate position for the sequence has positive numbers on a 
pair of adjacent nodes. 
A position $\lambda$ is {\em adjacency-free} if every game sequence 
played from $\lambda$ is adjacency-free.\myfootnote{For a dominant 
position $\lambda$, there can be both adjacency-free and 
non-adjacency-free game sequences.  For example, 
for the E-Coxeter graph\ \BThreeGraph\  
from the $\mathcal{B}_{3}$ family, the game sequence 
$(\gamma_{2}, \gamma_{1}, \gamma_{3}, \gamma_{2}, 
\gamma_{3}, \gamma_{1}, \gamma_{2})$ played from the fundamental 
position $\omega_{2}$ is adjacency-free   
while the game sequence 
$(\gamma_{2}, \gamma_{3}, \gamma_{2}, \gamma_{1}, 
\gamma_{2}, \gamma_{3}, \gamma_{2})$ is not adjacency-free.  Then the 
position $\omega_{2}$ for this E-Coxeter graph is not adjacency-free.}  
Following \S 1.1 of \cite{StemFC} and \S 8.1 of \cite{HumCoxeter},   
we let $\mathcal{W} = I_{n}^{*}$ 
be the free monoid on the set $I_{n}$.  
Elements of $\mathcal{W}$ are {\em words} and will be viewed as 
finite sequences of elements from $I_{n}$. The binary operation is 
concatenation, and the identity is the 
empty word. Fix a word 
$\mathbf{s} := ({i_{1}},\ldots,{i_{r}})$. 
Then $\myl_{\mathcal{W}}(\mathbf{s}) := r$ is the {\em length} of 
$\selt$.  
A {\em subword} of  
$\mathbf{s}$ 
is any subsequence 
$(i_{p},i_{p+1},\ldots,i_{q})$ of consecutive elements of $\selt$.  
For a nonnegative integer $m$ and $x,y \in I_{n}$, 
let $\langle x,y \rangle_{m}$ denote the sequence 
$(x,y,x,y,\ldots) \in \mathcal{W}$ so that 
$\myl_{\mathcal{W}}(\langle x,y \rangle_{m}) = m$. 
We employ several types of ``elementary simplifications'' in $\mathcal{W}$.  
An {\em elementary simplification of 
braid type} replaces a subword $\langle x,y \rangle_{m_{xy}}$ with 
the subword $\langle y,x \rangle_{m_{xy}}$ if $2 \leq m_{xy} < \infty$.   
An {\em elementary simplification of length-reducing type} 
replaces a subword $(x,x)$ with 
the empty subword. 
We let $\mathcal{S}(\mathbf{s})$ be the set of all words that can be 
obtained from $\mathbf{s}$ by some sequence of elementary 
simplifications of braid or length-reducing type.  
Since $s_{i}$ in $W$ is its own inverse for each $i \in I_{n}$, there is an 
induced mapping $\mathcal{W} \rightarrow W$.  We compose this with 
the mapping $W \rightarrow W$ for which $w \mapsto w^{-1}$ to get 
$\psi: \mathcal{W} \rightarrow W$ given by $\psi(\mathbf{s}) 
= s_{i_{r}}\cdots{s}_{i_{1}}$.  
Tits' Theorem for the word problem on Coxeter 
groups (cf.\ Theorem 8.1 of \cite{HumCoxeter}) 
implies that: {\sl For words $\mathbf{s}$ and $\mathbf{t}$ in 
$\mathcal{W}$, $\psi(\mathbf{s}) = \psi(\mathbf{t})$ if and only if 
$\mathcal{S}(\mathbf{s}) \cap \mathcal{S}(\mathbf{t}) \not= 
\emptyset$.} (This theorem is the basis for 
Part (1) of Eriksson's Reduced Word Result.) 
We say $\mathbf{s}$ is a {\em reduced} word for $w = \psi(\mathbf{s})$ if 
$\myl_{\mathcal{W}}(\mathbf{s}) = \myl(w)$ (assume this is the case 
for the remainder of the paragraph). Let $\mathcal{R}(w) 
\subseteq \mathcal{W}$ denote  
the set of all reduced words for $w$.  Suppose that $\mathbf{t} \in 
\mathcal{R}(w)$.  By Tits' Theorem, $\mathcal{S}(\mathbf{s}) \cap 
\mathcal{S}(\mathbf{t}) \not= \emptyset$, 
so that $\mathbf{t}$ can be obtained from 
$\mathbf{s}$ by a sequence of elementary simplifications of braid or 
length-reducing type.  Since $\myl_{\mathcal{W}}(\mathbf{s}) = \myl(w) = 
\myl_{\mathcal{W}}(\mathbf{t})$, then no elementary simplifications 
of length-reducing type can be used to obtain $\mathbf{t}$ from 
$\mathbf{s}$.  Then any member of $\mathcal{R}(w)$ can be obtained 
from any other member by a sequence of elementary simplifications of 
braid type.  
An {\em elementary simplification of commuting type} 
replaces a subword $(x,y)$ 
with the subword 
$(y,x)$ 
if $m_{xy} = 2$. 
The {\em commutativity class} $\mathcal{C}(\mathbf{s})$ of 
the word $\mathbf{s}$ is the set of all words that can be obtained 
from $\mathbf{s}$ by a sequence of elementary simplifications of 
commuting type.  
Clearly 
$\mathcal{C}(\mathbf{s}) \subseteq \mathcal{R}(w)$.  In fact there is a 
decomposition of $\mathcal{R}(w)$ into commutativity classes: 
$\mathcal{R}(w) = \mathcal{C}_{1} 
\disjointunion 
\cdots 
\disjointunion 
\mathcal{C}_{k}$, a disjoint union.  If 
$\mathcal{R}(w)$ has just one commutativity class, then $w$ is {\em 
fully commutative}.  Proposition 1.1 
of \cite{StemFC} states: {\sl An element $w \in W$ is fully 
commutative if and only if for all $x, y \in I_{n}$ such that $3 \leq 
m_{xy} < \infty$, there is no member of $\mathcal{R}(w)$ that 
contains $\langle x,y \rangle_{m_{xy}}$ as a subword}. 

\noindent 
{\bf \AdjacencyFreeForFinite}\ \ 
{\sl Let $J \subseteq I_{n}$. 
(1) Suppose $W_{J}$ is finite.  
Suppose an adjacency-free position $\lambda$ is $J^{c}$-dominant.  Then 
every element 
of $W^{J}$ is fully commutative. 
(2) Suppose $W$ is finite.  Suppose each element of $W^{J}$ is 
fully commutative.  Then any $J^{c}$-dominant position is adjacency-free.}

{\em Proof.} 
Our proof of (1) is by induction on the lengths of elements in 
$W^{J}$.  It is clear that the identity element is fully commutative. 
Now suppose that for all $v^{^{J}}$ in $W^{J}$ with $\myl(v^{^{J}}) < 
k$, it is the case that $v^{^{J}}$ is fully commutative, and 
consider $w^{^{J}}$ in $W^{J}$ such that 
$\myl(w^{^{J}}) = k$.  Suppose that for some adjacent $\gamma_{x}$ 
and $\gamma_{y}$ in $\Gamma$ with $3 \leq m_{xy} < \infty$, we have 
$\langle x,y \rangle_{m_{xy}}$ as a 
subword of some reduced word $\mathbf{s} = (i_{1},\ldots,i_{k}) \in 
\mathcal{R}(w^{^{J}})$.  Since $(i_{1},\ldots,i_{k-1})$ is a reduced 
word and $s_{i_{k-1}}\cdots{s}_{i_{1}}$ is in $W^{J}$, then 
$\langle x,y \rangle_{m_{xy}}$ cannot be a subword of 
$(i_{1},\ldots,i_{k-1})$.  
Therefore 
it must be 
the case that $\mathbf{s} = 
(i_{1},\ldots,i_{p},\langle x,y \rangle_{m_{xy}})$ for 
$p = k-m_{xy}$.  
Then, $\mathbf{s}' = (i_{1},\ldots,i_{p},\langle y,x \rangle_{m_{xy}})$ 
is also a reduced word for $w^{^{J}}$.  Since both $\mathbf{s}$ and 
$\mathbf{s}'$ correspond to legal firing sequences from $\lambda$ 
(\KeyWJResult), it 
must be the case that there are positive numbers at adjacent 
nodes $\gamma_{x}$ and $\gamma_{y}$ after the first $p$ firings.  But 
this contradicts the hypothesis that $\lambda$ is adjacency-free.  
Hence no reduced word for $w^{^{J}}$ can have a subword of the form 
$\langle x,y \rangle_{m_{xy}}$ for nodes 
$\gamma_{x}$ and $\gamma_{y}$ with $3 \leq m_{xy} < \infty$.  By 
Proposition 1.1 of \cite{StemFC} it follows that $w^{^{J}}$ is fully 
commutative, which completes the proof of part (1).

For part (2), assume every member of $W^{J}$ is fully commutative, and 
let $\lambda$ be any $J^{c}$-dominant position.  
Let $L := \myl(w_{0}) - 
\myl(\, (w_{0})_{_{J}}) = \myl(\, (w_{0})^{^{J}})$.  
Suppose an intermediate 
position $s_{i_{k}}\cdots{s}_{i_{1}}.\lambda$ 
for some game sequence $(\gamma_{i_{1}},\ldots,\gamma_{i_{L}})$ 
has positive numbers 
on adjacent nodes $\gamma_{x}$ and $\gamma_{y}$.  Then by Eriksson's 
Strong Convergence Theorem, there is a game sequence of length $L$ 
from $\lambda$ corresponding to a reduced word $\selt = 
(i_{1},\ldots,i_{k},\langle x,y 
\rangle_{m_{xy}},j_{k+m_{xy}+1},\ldots,j_{L})$ for $u := \psi(\selt)$.  
By \KeyWJResultCorollary, $u = (w_{0})^{^{J}}$. 
So $(w_{0})^{^{J}}$ is fully commutative (by 
hypothesis) and has reduced word $\selt$, in violation of Proposition 
1.1 of \cite{StemFC}.  Therefore $\lambda$ must be adjacency-free.\hfill\QED

In Theorem 5.1 of \cite{StemFC}, Stembridge classifies those $W^{J}$ 
for irreducible Coxeter groups $W$ such that every member of $W^{J}$ 
is fully commutative.  In view of \AdjacencyFreeForFinite\ and the 
classification of finite Coxeter groups, 
we may apply this result here to conclude that when 
$W$ is finite and irreducible, 
then the adjacency-free dominant positions of $(\Gamma,M)$ are exactly those 
specified in the following theorem. Observe that a dominant 
position $\lambda$ is adjacency-free if and only if $r\lambda := 
(r\lambda_{i})_{i \in I_{n}}$ is adjacency-free for all positive real 
numbers $r$. Call any such $r\lambda$ a {\em positive multiple} of 
$\lambda$.  

\noindent 
{\bf \AdjacencyFreeClassification}\ \ {\sl Suppose 
$(\Gamma,M)$ is connected.  If $W$ is 
finite, then an adjacency-free dominant position is a positive 
multiple of a fundamental position.  
All fundamental positions for any E-Coxeter graph of type 
$\mathcal{A}_{n}$ are adjacency-free.  The 
adjacency-free 
fundamental positions for any graph of type 
$\mathcal{B}_{n}$, 
$\mathcal{D}_{n}$, or $\mathcal{I}_{2}(m)$ are 
precisely those corresponding to end nodes. The adjacency-free 
fundamental positions for any graph of type 
$\mathcal{E}_{6}$, $\mathcal{E}_{7}$, 
or $\mathcal{H}_{3}$ are precisely those 
corresponding to the nodes marked with asterisks in 
\ECoxeterGraphFigure.  Any graph of type $\mathcal{E}_{8}$,  
$\mathcal{F}_{4}$, or $\mathcal{H}_{4}$ has no 
adjacency-free fundamental positions.}\hfill\QED 

For finite irreducible Coxeter groups $W$, it is a consequence of 
Theorems 5.1 and 6.1 of \cite{StemFC} that 
the Bruhat order $(W^{J},\leq)$ is a lattice 
if and only if $(W^{J},\leq)$ is a distributive lattice if and only if 
each element of $W^{J}$ is fully commutative. In these cases $|J^{c}| = 
1$ and all such $J^{c}\, '${\small s} correspond to the adjacency-free fundamental 
positions from \AdjacencyFreeClassification\ above.  
\AdjacencyFreeForFinite\ above adds to 
these equivalences the property that each element of $W^{J}$ is fully 
commutative if and only if for any associated E-GCM graph, 
any $J^{c}$-dominant position is adjacency-free.  The adjacency-free 
viewpoint is similar to Proctor's original viewpoint (cf.\ Lemma 3.2 
of \cite{PrEur}).

%==================================================================
\vspace{1ex} 

\noindent
{\Large \bf \RootsNum.\ \ 
Generating positive roots from E-game play}

\vspace{1ex} 
The results of this section 
expand on Remark 4.6 of \cite{ErikssonDiscrete}. 
The goal here is to characterize when all positive roots can be 
obtained from a single game sequence, 
as in the following example.  In \BtwoMarsDemo\ 
with amplitude matrix $M = \left(\begin{array}{cc} 2 & -1\\ -2 & 
2\end{array}\right)$, assume the initial position $\lambda = (a,b)$ is 
strongly dominant.  For the game sequence $(\gamma_{2}, \gamma_{1}, 
\gamma_{2}, \gamma_{1})$, notice that the respective numbers at 
the fired nodes are $b$, $a+2b$, $a+b$, and $a$.  Thought of now as root 
functionals, the latter  
are in one-to-one correspondence with the 
positive roots $\Phi_{M}^{+} = \{\alpha_{2}, \alpha_{1}+2\alpha_{2}, 
\alpha_{1}+\alpha_{2}, \alpha_{1}\}$.  For  
$M = \left(\begin{array}{cc} 2 & -1/2\\ -2 & 
2\end{array}\right)$ with E-Coxeter graph in the $\mathcal{A}_{2}$ 
family (cf.\ Exercise 4.9 of \cite{BB}), the situation is different.  
From a strongly dominant position $\lambda = (a,b)$ on \WeirdATwoEGCMGraphh, the game 
sequence $(\gamma_{2}, \gamma_{1}, \gamma_{2})$ has respective 
numbers $b$, $a+2b$, and $\frac{1}{2}a$ at the fired nodes.  
However, the positive roots are 
$\Phi_{M}^{+} = \{\alpha_{2}, \alpha_{1}+2\alpha_{2}, 
\frac{1}{2}\alpha_{1}, \alpha_{1}, \frac{1}{2}\alpha_{1}+\alpha_{2}, 
2\alpha_{2}\}$. 

In general, for $p \geq 1$  
suppose $\mathbf{s} := (\gamma_{i_{1}},\ldots,\gamma_{i_{p}})$ is a legal 
firing sequence from some initial position $\lambda$ on $(\Gamma,M)$.  After 
$(\gamma_{i_{1}},\ldots,\gamma_{i_{q-1}})$ is played ($1 \leq q \leq 
p$), 
the number at node $\gamma_{i_{q}}$ is 
$\langle s_{i_{q-1}}\cdots{s_{i_{1}}}.\lambda, \alpha_{i_{q}} \rangle = 
\langle \lambda, s_{i_{1}}\cdots{s_{i_{q-1}}}.\alpha_{i_{q}} \rangle = 
\phi_{\beta_{q}}(\lambda)$ with $\beta_{q} := 
s_{i_{1}}\cdots{s_{i_{q-1}}}.\alpha_{i_{q}}$.  
With $\mathbf{s}$ and $\lambda$ understood, then we  
say $\phi_{\beta_{q}}$ is the {\em root functional at node} 
$\gamma_{i_{q}}$.\myfootnote{It follows from Part (2) of Eriksson's 
Reduced Word Result and \EquivalenceRemark\ that for 
any given strongly dominant position $\lambda$ 
and any 
positive root $\alpha$, there is a legal firing sequence 
$(\gamma_{i_{1}},\ldots,\gamma_{i_{q-1}})$ 
played from $\lambda$ 
such that $\phi_{\alpha}$ is the root functional at node 
$\gamma_{i_{q}}$.} 
By Part (1) of 
Eriksson's Reduced Word Result, $w := s_{i_{p}}\cdots{s_{i_{2}}}s_{i_{1}}$ is 
reduced. 
This is exactly 
the situation of Exercise 5.6.1 of \cite{HumCoxeter}, where the 
representation is the ``standard'' geometric representation of $W$. 
There, one concludes that the 
$\beta_{q}\, '${\small s} are 
distinct and precisely all of the positive roots $\beta$ for which 
$w.\beta$ is a negative root. 
In our more general setting we have: 

\noindent 
{\bf \NoRepeatPosRootsProp}\ \ 
{\sl Let $w = s_{i_{p}}\cdots{s_{i_{2}}}s_{i_{1}}$ with $\myl(w) = p \geq 1$.  Let 
$\beta_{q} := s_{i_{1}}s_{i_{2}}\cdots{s_{i_{q-1}}}.\alpha_{i_{q}}$ 
for $1 \leq q \leq p$.  
Then $\beta_{q} \not= \beta_{r}$ for $q \not= r$ and 
$\{\beta_{q}\}_{q=1}^{p} \subseteq N_{M}(w)$.  Moreover, 
$\{\beta_{q}\}_{q=1}^{p} = N_{M}(w)$ if and only if for $1 \leq q \leq 
p$ the ON-connected component $(\Gamma',M')$ containing 
$\gamma_{i_{q}}$ is unital ON-cyclic with $f_{\Gamma',M'} = 1$.} 

{\em Proof.} 
Each $\beta_{q} \in \Phi_{M}^{+}$ by \EquivalenceRemark. 
Also, 
$w.\beta_{q} = s_{i_{p}}s_{i_{p-1}}\cdots{s_{i_{q}}}.\alpha_{i_{q}} 
\in \Phi_{M}^{-}$ follows from the fact that  
$\myl(s_{i_{p}}s_{i_{p-1}}\cdots{s_{i_{q+1}}}s_{i_{q}}s_{i_{q}}) <  
\myl(s_{i_{p}}s_{i_{p-1}}\cdots{s_{i_{q+1}}}s_{i_{q}})$. Hence 
$\beta_{q} \in N_{M}(w)$. 
For $q < r$, suppose $\beta_{q} = \beta_{r}$.  
Then one can see that 
$s_{i_{q}}\cdots{s}_{i_{r-1}}.\alpha_{i_{r}} = \alpha_{i_{q}}$, and 
so $s_{i_{q+1}}\cdots{s}_{i_{r-1}}.\alpha_{i_{r}} = -\alpha_{i_{q}} 
\in \Phi_{M}^{-}$. Then $\myl(s_{i_{q+1}}\cdots{s}_{i_{r-1}}s_{i_{r}}) 
< \myl(s_{i_{q+1}}\cdots{s}_{i_{r-1}})$.  
But $s_{i_{q+1}}\cdots{s}_{i_{r-1}}s_{i_{r}}$ is reduced 
and longer than $s_{i_{q+1}}\cdots{s}_{i_{r-1}}$, a contradiction.   
So $\beta_{q} \not= \beta_{r}$.
For the ``if'' direction of the last assertion of the lemma, by 
\LengthProposition\ $N_{M}(w)$ is finite.  Since $f_{1} = f_{2} = 1$, 
then $\myl(w) = |N_{M}(w)| = p$.  For the ``only if'' direction, 
$N_{M}(w)$ has finite order $p = \myl(w)$.  Then by \LengthProposition, each 
$\mathfrak{S}_{M}(\alpha_{i_{q}})$ is finite, so by \TFAE\ the 
ON-connected component $(\Gamma',M')$ containing $\gamma_{i_{q}}$ is 
unital ON-cyclic.  Combining $\myl(w) = |N_{M}(w)|$ and 
$f_{1}\myl(w) \leq |N_{M}(w)| \leq f_{2}\myl(w)$ gives $f_{1} = 
f_{2} = 1$.  Therefore $f_{\Gamma',M'} = 1$.\hfill\QED

From this lemma, it is apparent now why the game sequence exhibited 
in the above $\mathcal{A}_{2}$ example failed to generate all of the 
positive roots: the E-GCM graph has an odd asymmetry which results 
in some positive roots which are nontrivial multiples of simple 
roots.  In this case, $f_{\Gamma,M} = 2 = 
|\mathfrak{S}_{M}(\alpha_{i})|$ for $i = 1,2$.  
The positive roots $\{\alpha_{2}, \alpha_{1}+2\alpha_{2}, 
\frac{1}{2}\alpha_{1}\}$ associated with the root functionals of 
the game sequence 
$(\gamma_{2}, \gamma_{1}, \gamma_{2})$ are a proper subset of 
$N_{M}(w_{0}) = \Phi_{M}^{+}$ where $w_{0} = s_{2}s_{1}s_{2}$.  
In general, if $W$ is finite and $(\Gamma,M)$ has odd asymmetries  
then not every positive root will be encountered 
as a positive root functional in a given game sequence, as the 
next result shows. 
However, if the amplitude matrix $M$ is integral, then $(\Gamma,M)$ 
has no odd asymmetries and thus enjoys the equivalent 
properties of the following theorem. 

\noindent 
{\bf \AllPositiveRootsProp}\ \ {\sl Suppose $W$ is finite.  Let 
$s_{i_{l}}\cdots{s_{i_{2}}}s_{i_{1}}$ be any reduced expression for 
$w_{0}$. 
For $1 \leq j \leq l$, set $\beta_{j} := 
s_{i_{1}}s_{i_{2}}\cdots{s_{i_{j-1}}}.\alpha_{i_{j}}$.  
Then the following are equivalent:} 

\hspace{0.05in} {\sl (1) $(\Gamma,M)$ 
has no odd asymmetries;} 
 
\hspace{0.05in} {\sl (2) Each  
ON-connected component $(\Gamma',M')$ of $(\Gamma,M)$ is unital 
ON-cyclic with $f_{\Gamma',M'} = 1$;} 

\hspace{0.05in} {\sl (3) $\{\beta_{j}\}_{j=1}^{l} = \Phi_{M}^{+}$;}

\hspace{0.05in} {\sl (4) $\myl(w_{0}) = |\Phi_{M}^{+}|$;} 

\hspace{0.05in} {\sl (5) Each positive root appears 
as the root functional 
$\phi_{\beta_{j}}$ at some node $\gamma_{i_{j}}$ for the game  
sequence $(\gamma_{i_{1}},\ldots,\gamma_{i_{l}})$ played from any 
strongly dominant position.}

{\em Proof.} 
For (1) $\Leftrightarrow$ (2), 
note that by \HumphreysTheorem\ a nontrivial positive multiple 
of some simple root is itself a root if and only if there are odd asymmetries.  
For (2) $\Rightarrow$ (3), recall from \S \ConvergenceNum\ that
$N_{M}(w_{0}) = \Phi_{M}^{+}$. \NoRepeatPosRootsProp\ shows that 
that $\{\beta_{j}\}_{j=1}^{l} = N_{M}(w_{0})$, so (3) follows. 
(4) follows immediately from (3).  For (4) 
$\Rightarrow$ (5), first note that by Part (2) of Eriksson's Reduced 
Word Result, the firing sequence $(\gamma_{i_{1}},\ldots,\gamma_{i_{l}})$ 
is legal from any strongly dominant position, and by 
\LengthUsingLongestWord\ this is a game sequence.  \NoRepeatPosRootsProp\ and 
comments preceding that lemma show that for this game 
sequence 
the positive roots in the set $\{\beta_{j}\}_{j=1}^{l}$ appear 
precisely once each as root functionals.  The hypothesis  
$\myl(w_{0}) = |\Phi_{M}^{+}|$ means that  $\{\beta_{j}\}_{j=1}^{l} = 
\Phi_{M}^{+}$, from which (5) follows.  
To show (5) $\Rightarrow$ (2), 
choose an ON-connected component  $(\Gamma',M')$.  \PropList\ 
show that $(\Gamma',M')$ must be unital ON-cyclic, else $W$ will 
be infinite. 
Let $J$ be the subset of 
$I_{n}$ corresponding to the nodes of the subgraph $\Gamma'$.  
For notational 
convenience set $w = w_{0}$, $w_{_{J}} = (w_{0})_{_{J}}$, and 
$w^{^{J}} = (w_{0})^{^{J}}$.  
Set $w_{_{J}} = s_{j_{k}}\cdots{s}_{j_{2}}s_{j_{1}}$, a reduced expression.  
Using \PositiveToNegativeRootsLemma, we see that 
\begin{eqnarray*} 
|N_{M}(ws_{j_{1}})| & = & |N_{M}(w)| - f_{\Gamma',M'},\\
|N_{M}(ws_{j_{1}}s_{j_{2}})| & = & |N_{M}(ws_{j_{1}})| - 
f_{\Gamma',M'}\ \, =\ \, |N_{M}(w)| - 2f_{\Gamma',M'},
\end{eqnarray*}
so that eventually $|N_{M}(w)| = |N_{M}(w^{^{J}})| + 
\myl(w_{_{J}})f_{\Gamma',M'}$.  Now by hypothesis 
each positive root functional 
appears once and therefore, by \NoRepeatPosRootsProp, exactly once.  
Then $l = \myl(w) = |\Phi_{M}^{+}| = |N_{M}(w)|$. 
By \LengthProposition, $|N_{M}(w^{^{J}})| \geq \myl(w^{^{J}})$.  
Summarizing, 
$\myl(w^{^{J}}) + \myl(w_{_{J}}) = \myl(w) = |N_{M}(w)| = 
|N_{M}(w^{^{J}})| + \myl(w_{_{J}})f_{\Gamma',M'} \geq 
\myl(w^{^{J}}) + \myl(w_{_{J}})f_{\Gamma',M'}$,
from which $f_{\Gamma',M'} = 1$.  
\hfill\QED

%==================================================================
\vspace{1ex} 

\noindent
{\Large \bf \DynkinClassNum.\ \ A Dynkin diagram classification of 
E-GCM graphs meeting a certain finiteness requirement}

\vspace{1ex} 
%===========================
% Non-unital-ON-cyclic stuff
%===========================
We say a connected E-GCM graph is {\em admissible} if 
there exists a nonzero dominant initial position with a convergent 
game sequence.  In this section we prove the following Dynkin diagram 
classification result. 

\noindent 
{\bf \EGCMTheorem}\ \ {\sl A connected E-GCM graph is 
admissible if and only if it 
is a connected E-Coxeter graph.  
In these cases, for any given initial position every game sequence  
will converge to the same terminal position in the 
same finite number of steps.} 

Our proof of \EGCMTheorem\ given at the end of this section 
uses the classification of finite 
Coxeter groups.  Another proof based on ideas from 
\cite{ErikssonThesis} is given in \cite{DE}. 
That proof uses  
combinatorial reasoning together with a 
result from the Perron--Frobenius 
theory for eigenvalues of nonnegative real matrices, and it  
does not require the classification of 
finite Coxeter groups. 
Before proceeding toward our proof of \EGCMTheorem, we record 
two closely related results. 
In \cite{ErikssonThesis}, Eriksson establishes 
the following result.  
(For an ``A-D-E'' version, see \cite{ErikssonLinear}.) 
The statement we give here essentially 
combines his Theorems 6.5 and 6.7.  
An E-GCM graph is {\em strongly admissible} if every nonzero 
dominant position has a convergent game sequence.  

\noindent 
{\bf \ErikssonTheorem\ (Eriksson)}\ \ {\sl A connected E-GCM graph 
is strongly admissible if and only if it is a connected E-Coxeter 
graph.}

Using this result Eriksson re-derives in \S 8.4 of 
\cite{ErikssonThesis} the well-known classification of finite 
irreducible Coxeter groups, which we state as: {\sl 
An irreducible Coxeter group 
$W(\Gamma,M)$ is finite if and only if the connected E-GCM graph 
$(\Gamma,M)$ is an E-Coxeter graph from \ECoxeterGraphFigure.}  
In Propositions 4.1 and 4.2 of \cite{Deodhar}, Deodhar gives a number of 
statements equivalent to the assertion that a given irreducible 
Coxeter group is finite.  
As an immediate consequence of \TwoTheorems\ and the classification of 
finite irreducible Coxeter groups, we add to that list 
the following equivalence.  

\noindent 
{\bf \DeodharEquivalenceCorollary}\ \ 
{\sl An irreducible Coxeter group $W$ is finite if and only if 
there is an admissible E-GCM graph whose associated Coxeter group is 
$W$ if and only if any E-GCM graph is strongly admissible 
when its associated Coxeter group is $W$.}\hfill\QED 

Extending 
\KeyWJResult\ to all subsets 
$J \subseteq I_{n}$ would yield a simple proof of the first assertion 
of \EGCMTheorem: For any 
given proper subset $J \subset I_{n}$, the E-GCM 
graph $(\Gamma,M)$ would have a convergent game sequence for some 
$J^{c}$-dominant $\lambda$ if and only if $W^{J}$ is finite if and only if 
$W$ is finite (by Proposition 4.2 of \cite{Deodhar}).  
Observe that the ``if'' direction of the first 
assertion in \EGCMTheorem\ follows from \ErikssonTheorem.  
The second assertion in \EGCMTheorem\ follows from Eriksson's 
Strong Convergence Theorem.  
So our effort in the proof of \EGCMTheorem\ will be 
mainly concerned with  
demonstrating the ``only if'' part of the first assertion.   
Our proof of this part 
is by induction on the number of nodes.  
The main idea of our proof is to use reductions effected by the 
preliminary results of Section \PrelimNum\ together with some further 
results derived here.  The lemmas that follow use 
\NotMarsFriendlyLemma, which depends crucially on Eriksson's 
Comparison Theorem. 
We say an $n$-node graph $\Gamma$ is a {\em loop} if the nodes 
can be numbered $\gamma_{1},\ldots,\gamma_{n}$  in such a way that 
for all $1 \leq i \leq n$, 
$\gamma_{i}$ is adjacent precisely to $\gamma_{i+1}$ 
and $\gamma_{i-1}$, understanding that $\gamma_{0} = \gamma_{n}$ 
and $\gamma_{n+1} = \gamma_{1}$. 

\noindent 
{\bf \CycleLemma}\ \ 
{\sl Suppose that the underlying 
graph $\Gamma$ of an E-GCM graph $(\Gamma,M)$ is a loop 
and that for any edge in $(\Gamma,M)$ 
the amplitude product is unity.  
Then $(\Gamma,M)$ is not admissible.} 

{\em Proof.} We find a divergent 
game sequence starting from 
the fundamental position $\omega_{1}$. Then by renumbering the nodes, 
we see that every fundamental position will have a divergent game 
sequence, and by \NotMarsFriendlyLemma\ it then follows that 
$(\Gamma,M)$ is not admissible. 
Let the ON-cycle $\mathcal{C}$ 
be $[\gamma_{1},\gamma_{2},\ldots,\gamma_{n},\gamma_{1}]$. 
From initial position $\omega_{1}$ 
we propose starting with the firing 
sequence 
$(\gamma_{1},\ldots,\gamma_{n-1},
\gamma_{n},\gamma_{n-1},\ldots\gamma_{2})$.  One can check that all 
of these node firings are legal and that the resulting numbers are zero 
at all nodes 
other than $\gamma_{1}$, $\gamma_{2}$, and $\gamma_{n}$.  The 
numbers at the latter nodes are, respectively, 
$1+\Pi_{_{\mathcal{C}}}+\Pi_{_{\mathcal{C}}}^{-1}$, 
$M_{12}(\Pi_{_{\mathcal{C}}}^{-1})$, and 
$M_{1n}(\Pi_{_{\mathcal{C}}})$. 
By repeating the proposed firing sequence 
$(\gamma_{1},\ldots,\gamma_{n-1},
\gamma_{n},\gamma_{n-1},\ldots\gamma_{2})$ from this position 
we obtain zero at all nodes except at $\gamma_{1}$, $\gamma_{2}$, and 
$\gamma_{n}$, which are now 
$1+\Pi_{_{\mathcal{C}}}+\Pi_{_{\mathcal{C}}}^{-1}+
\Pi_{_{\mathcal{C}}}^{2}+\Pi_{_{\mathcal{C}}}^{-2}$, 
$M_{12}(\Pi_{_{\mathcal{C}}}^{-1}+\Pi_{_{\mathcal{C}}}^{-2})$, and 
$M_{1n}(\Pi_{_{\mathcal{C}}}+\Pi_{_{\mathcal{C}}}^{2})$ 
respectively.  After $k$ applications of the proposed firing sequence 
we have numbers 
$1+\sum_{j=1}^{k}\Pi_{_{\mathcal{C}}}^{j}+\Pi_{_{\mathcal{C}}}^{-j}$, 
$M_{12}(\sum_{j=1}^{k}\Pi_{_{\mathcal{C}}}^{-j})$, and 
$M_{1n}(\sum_{j=1}^{k}\Pi_{_{\mathcal{C}}}^{j})$ 
at nodes 
$\gamma_{1}$, $\gamma_{2}$, and 
$\gamma_{n}$, and zeros elsewhere.  Thus we have exhibited 
a divergent game sequence.\hfill\QED

\noindent 
{\bf \FourCycleLemma}\ \ 
{\sl An E-GCM graph in the family} 
\hspace*{0.1in}
\parbox[c]{0.75in}{
\setlength{\unitlength}{0.5in} 
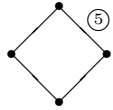
\begin{picture}(1,1.35)
            \put(0,0.6){\circle*{0.075}}
            \put(0.5,0.1){\circle*{0.075}}
            \put(0.5,1.1){\circle*{0.075}}
            \put(1,0.6){\circle*{0.075}}
            %=======
            \put(0,0.6){\line(1,1){0.5}}
            \put(0,0.6){\line(1,-1){0.5}}
            \put(1,0.6){\line(-1,1){0.5}}
            \put(1,0.6){\line(-1,-1){0.5}}
            %=======
            \put(0.6,0.85){\CircleInteger{5}}
           \end{picture}
}{\sl is not admissible.}\hfill\QED 

{\em Proof.} Let $(\Gamma,M)$ be an E-GCM graph in the given family.  
Label the nodes $\gamma_{1}$, $\gamma_{2}$, 
$\gamma_{3}$, and $\gamma_{4}$ clockwise from the top.  
Our strategy is to show that the repeating firing sequence 
$\mathbf{r} := (\mathbf{s}, \mathbf{s}, \ldots)$ 
can be legally applied to some position 
obtained from E-game play starting with any given fundamental 
position, where $\mathbf{s}$ is the 
subsequence  
$(\gamma_{1}, \gamma_{2}, 
\gamma_{3}, \gamma_{4})$. 
This will give us a divergent game sequence from each fundamental 
position, so by \NotMarsFriendlyLemma\ it will follow that 
$(\Gamma,M)$ is not admissible. 
For adjacent 
nodes $\gamma_{1}$ and $\gamma_{2}$, set $p := -M_{12}$, 
$q := -M_{21}$.  Note that $pq = (3+\sqrt{5})/2$.  
Set $r := -M_{23}$, $s := -M_{32}$, $t := -M_{34}$, $u := -M_{43}$, 
$v := -M_{41}$, and $w := -M_{14}$.  
We have $rs = tu = vw = 1$.  
Note that $p, q, r, s, t, u, v, 
w$ are the absolute values of the amplitudes read in alphabetical 
order clockwise from the top. 
We say a position 
$(a,b,c,d)$ meets condition ({\tt *}) if $a > 0$, $b \geq 0$, $c \geq 
0$, $d \leq 0$, $aw+d \geq 0$, and $aprt + brt + ct + d > 0$.  One 
can easily check that from any such position 
the firing sequence $\mathbf{s} = 
(\gamma_{1}, \gamma_{2}, 
\gamma_{3}, \gamma_{4})$ is legal: The positive numbers at the 
fired nodes are respectively $a$, $ap+b$, $apr+br+c$, and 
$aw+aprt+brt+ct+d$.  The resulting position is 
$(A,B,C,D)$ 
with $A = \frac{3+\sqrt{5}}{2}a + bq + v(aprt + brt + ct + d)$, 
$B = sc$, $C = u(aw+d)$, and $D = -aw - aprt - brt - ct - d$.  
Clearly $A > 0$, $B \geq 0$, $C \geq 0$, and $D < 0$.  Also, $Aw + 
D = (\frac{3+\sqrt{5}}{2} - 1)aw + bqw > 0$, and $Aprt + Brt + Ct + D = 
(\frac{3+\sqrt{5}}{2} - 1)aprt + (\frac{3+\sqrt{5}}{2} - 1)brt + 
prtv(aprt + brt + ct + d) > 0$.  So, $(A,B,C,D)$ meets condition 
({\tt *}).  The fundamental position $\omega_{1} = (1,0,0,0)$ 
meets condition ({\tt *}), so it follows that the divergent firing 
sequence $\mathbf{r}$ can be legally played from 
this initial position.  Play the legal sequence $(\gamma_{2}, 
\gamma_{3}, \gamma_{4})$ from the fundamental position $\omega_{2} = 
(0,1,0,0)$ to obtain the position $(q+rtv,0,0,-rt)$.  It is easily 
checked that the latter position meets condition ({\tt *}).  It 
follows that the divergent firing 
sequence $(\gamma_{2}, \gamma_{3}, \gamma_{4}, \mathbf{r})$ can be 
legally played from $\omega_{2}$. 
Similarly see that the divergent firing 
sequence $(\gamma_{3}, \gamma_{4}, \mathbf{r})$ can be 
legally played from $\omega_{3}$ and that the divergent firing 
sequence $(\gamma_{4}, \mathbf{r})$ can be 
legally played from $\omega_{4}$.\hfill\QED

\noindent 
{\bf \ThreeCycleLemma}\ \ 
{\sl Suppose $(\Gamma,M)$ is the following three-node E-GCM graph:}  
\hspace*{0.1in}
\parbox[c]{0.5in}{
\setlength{\unitlength}{0.75in}
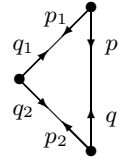
\begin{picture}(0.6,1.2)
            \put(0,0.6){\circle*{0.075}}
            \put(0.5,0.1){\circle*{0.075}}
            \put(0.5,1.1){\circle*{0.075}}
            %=======
            \put(0,0.6){\line(1,1){0.5}}
            \put(0,0.6){\line(1,-1){0.5}}
            \put(0.5,0.1){\line(0,1){1}}
            %=======
            \put(0.5,0.3){\vector(0,1){0.1}}
            \put(0.6,0.3){\footnotesize $q$}
            \put(0.5,0.9){\vector(0,-1){0.1}}
            \put(0.6,0.8){\footnotesize $p$}
            \put(0,0.6){\vector(1,1){0.2}}
            \put(0.5,1.1){\vector(-1,-1){0.2}}
            \put(0,0.6){\vector(1,-1){0.2}}
            \put(0.5,0.1){\vector(-1,1){0.2}}
            \put(-0.05,0.8){\footnotesize $q_{1}$}
            \put(0.175,1){\footnotesize $p_{1}$}
            \put(-0.05,0.35){\footnotesize $q_{2}$}
            \put(0.175,0.15){\footnotesize $p_{2}$} 
\end{picture}
} 
\hspace*{0.25in}{\sl Assume that all node pairs are odd-neighborly.  
Then $(\Gamma,M)$ is not admissible.} 

{\em Notes on the proof.} As in the proofs of \RuleOutSmallCyclesLemmas, 
we apply \NotMarsFriendlyLemma\ after showing that 
from each fundamental position there is a legal firing sequence that 
can be repeated indefinitely.  However, the variable amplitude 
products on edges of this graph make this argument a little more 
delicate than our arguments for the previous lemmas.  
A key part of the argument in 
this case is an explicit computation of matrix representations of 
powers of 
$\sigma_{M}(s_{i}s_{j})$ with respect to the basis 
$\{\alpha_{1}, \alpha_{2}, \alpha_{3}\}$ of simple roots. 
These computations are used to 
understand positions resulting from alternating sequences of 
firings on adjacent nodes. 
For complete details, see \cite{DonSupplement}.

We can now prove \EGCMTheorem. 

{\em Proof of \EGCMTheorem.}  
First we use induction on $n$, the number of nodes, to show that any 
connected admissible E-GCM graph must be from one of the 
families of \ECoxeterGraphFigure.  
Clearly a one-node E-GCM graph is admissible.  
For some $n \geq 2$, 
suppose the result is true for all connected admissible E-GCM 
graphs with fewer than $n$ nodes.  Let $(\Gamma,M)$ be a connected, 
admissible, $n$-node E-GCM graph.  
Suppose $(\Gamma,M)$ is 
unital ON-cyclic.  Then 
by \RuleOutInfinitePropositions, we must have $W$ finite.  
Then by the classification of finite irreducible Coxeter groups, 
$(\Gamma,M)$ must be in one of the families of graphs in  
\ECoxeterGraphFigure.  
Now suppose $(\Gamma,M)$ is not unital ON-cyclic.  
First we show that any cycle (ON or otherwise) in $(\Gamma,M)$ must 
use all $n$ nodes.  Indeed, the (connected) 
E-GCM subgraph $(\Gamma',M')$ whose nodes are the nodes of 
a cycle must be admissible by  
\SubgraphLemma.  If $(\Gamma',M')$ has fewer than $n$ nodes, then 
the induction hypothesis applies. But E-Coxeter graphs have no 
cycles (ON or otherwise), so $(\Gamma',M')$ must be all  
of $(\Gamma,M)$. 
Second, $(\Gamma,M)$ 
has an ON-cycle $\mathcal{C}$ for which $\Pi_{_{\mathcal{C}}} \not= 
1$. We can make the following choice for $\mathcal{C}$: Choose 
$\mathcal{C}$ to be a simple ON-cycle with $\Pi_{_{\mathcal{C}}} 
\not= 1$ whose length is as small as possible.  This smallest length 
must therefore be $n$.  
We wish to show that the underlying graph $\Gamma$ 
is a loop.  Let the numbering of the nodes of $\Gamma$ follow 
$\mathcal{C}$, so 
$\mathcal{C} = [\gamma_{1}, \gamma_{2}, \ldots, \gamma_{n}, 
\gamma_{1}]$.  If $\Gamma$ is not a loop, then there are adjacencies 
amongst the $\gamma_{i}$$'${\footnotesize s} besides those of 
consecutive elements of 
$\mathcal{C}$.  But this in turn means that $(\Gamma,M)$ has a 
cycle that uses fewer than $n$ nodes.   So $\Gamma$ is a loop.  Of 
course we must have $n \geq 3$. 
\ThreeCycleLemma\ rules out the possibility that $n=3$.  
Any connected  
E-GCM subgraph $(\Gamma',M')$ obtained from $(\Gamma,M)$ by removing 
a single node must now be a ``branchless'' E-Coxeter graph from 
\ECoxeterGraphFigure\ whose adjacencies are all odd.  So if $n=4$, 
$(\Gamma,M)$ must be in one of the families 
\hspace*{0.1in}
\parbox[c]{0.75in}{
\setlength{\unitlength}{0.5in} 
\begin{picture}(1,1.35)
            \put(0,0.6){\circle*{0.075}}
            \put(0.5,0.1){\circle*{0.075}}
            \put(0.5,1.1){\circle*{0.075}}
            \put(1,0.6){\circle*{0.075}}
            %=======
            \put(0,0.6){\line(1,1){0.5}}
            \put(0,0.6){\line(1,-1){0.5}}
            \put(1,0.6){\line(-1,1){0.5}}
            \put(1,0.6){\line(-1,-1){0.5}}
           \end{picture}
} 
or 
\hspace*{0.1in}
\parbox[c]{0.75in}{
\setlength{\unitlength}{0.5in} 
\begin{picture}(1,1.35)
            \put(0,0.6){\circle*{0.075}}
            \put(0.5,0.1){\circle*{0.075}}
            \put(0.5,1.1){\circle*{0.075}}
            \put(1,0.6){\circle*{0.075}}
            %=======
            \put(0,0.6){\line(1,1){0.5}}
            \put(0,0.6){\line(1,-1){0.5}}
            \put(1,0.6){\line(-1,1){0.5}}
            \put(1,0.6){\line(-1,-1){0.5}}
            %=======
            \put(0.6,0.85){\CircleInteger{5}}
           \end{picture}
}\hspace*{-0.1in}, which are ruled out by \RuleOutSmallCyclesLemmas\ 
respectively.  
If $n \geq 5$, the only possibility is that $(\Gamma,M)$ meets 
the hypotheses of \CycleLemma\ and therefore 
is not admissible.  In all cases, 
we see that if $(\Gamma,M)$ is not unital ON-cyclic, then it is 
not admissible.  This completes the induction step, so we have shown 
that a connected admissible E-GCM graph must be in one of the 
families of \ECoxeterGraphFigure. 

On the other hand, if $(\Gamma,M)$ is 
from \ECoxeterGraphFigure, then the Coxeter group $W$ 
is finite (again by the classification), 
so there is an upper bound on the length of any element in $W$. 
So by 
Part (1) of Eriksson's Reduced Word Result, the numbers game 
converges for any initial position. 
The remaining claims of \EGCMTheorem\ now follow from 
Eriksson's Strong Convergence Theorem.\hfill\QED

\vspace*{-0.15in}

%=============================================
% Bibliography
%=============================================
\renewcommand{\baselinestretch}{1}

\end{document}